\documentclass[10pt]{article}
\usepackage{amsfonts, amssymb, amsmath, amscd, longtable, tabularx, multicol}
\usepackage{graphicx}
\newtheorem{theorem-intro}{Theorem}[section]
\newtheorem{theorem}{Theorem}
\newtheorem{lemma}{Lemma}

\newtheorem{corollary}{Corollary}
\newtheorem{corollary-intro}[theorem-intro]{Corollary}
\newtheorem{remark}{Remark}

\newtheorem{conjecture}{Conjecture}
\newtheorem{example}{Example}
\newtheorem{question}{Question}
\def\Z{\mathbb Z}
\def\R{\mathbb R}

\def\pf {\noindent {\bf Proof:} \ }
\def\be{\begin{enumerate}}
\def\ee{\end{enumerate}}
\def\bi{\begin{itemize}}
\def\ei{\end{itemize}}

\long\def\symbolfootnote[#1]#2{\begingroup\def\thefootnote{\fnsymbol{footnote}}\footnote[#1]{#2}\endgroup}

\renewcommand{\c}{\cite}
\newcommand{\s}{\sigma_}
\newcommand{\done}{\hfill $\blacksquare$}

\newcommand{\TT}{\mathbb T}
\renewcommand{\L}{\mathbb L}
\newcommand{\T}{\mathcal T}
\newcommand{\kb}[1]{\ensuremath{\langle\, #1 \,\rangle}}

\begin{document}

\title{A new twist on Lorenz links}

\author{
Joan Birman
\footnote{Supported in part by National Science Foundation DMS-0405586.}\\
{\em {\small Barnard College, Columbia University}}\\ \\
 Ilya Kofman
\footnote{Supported by National Science Foundation DMS-0456227 and a PSC-CUNY grant.}\\
{\em {\small College of Staten Island, City University of New York}}}
\date{April 15, 2009}
\maketitle

\begin{abstract} \noindent 
  Twisted torus links are given by twisting a subset of strands on a
  closed braid representative of a torus link.  T--links are a natural
  generalization, given by repeated positive twisting.  We establish a
  one-to-one correspondence between positive braid representatives of
  Lorenz links and T--links, so Lorenz links and T--links coincide.
  Using this correspondence, we identify over half of the simplest
  hyperbolic knots as Lorenz knots.  We show that both hyperbolic
  volume and the Mahler measure of Jones polynomials are bounded for
  infinite collections of hyperbolic Lorenz links.  The correspondence
  provides unexpected symmetries for both Lorenz links and T-links,
  and establishes many new results for T-links, including new braid
  index formulas.
\end{abstract}
\symbolfootnote[0]{2000 Mathematics Subject Classification 57M25, 57M27, 57M50}

\section{Introduction}
\label{S:introduction}

The Lorenz differential equations \cite{Lo} have become well-known as
the prototypical chaotic dynamical system with a ``strange attractor''
(see \c{Viana}, and references therein).  A periodic orbit in the flow
on $\mathbb{R}^3$ determined by the Lorenz equations is a closed curve
in $\mathbb{R}^3$, which defines a {\em Lorenz knot}.  Lorenz knots
contain many known classes of knots, but the complete classification
of Lorenz knots remains open: What types of knots can occur?

\begin{figure}[htpb!]
\centerline{\includegraphics[scale=.7] {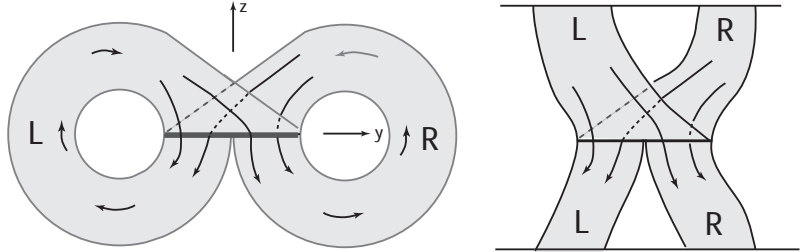}}
\caption{(a)The Lorenz template. \ (b)The Lorenz braid template}
\label{F:template}
\end{figure}

Guckenheimer and Williams introduced the {\em Lorenz template}, also
called the {\em geometric Lorenz attractor}, which is an embedded
branched surface in $\mathbb{R}^3$ with a semi-flow.  Later, Tucker
\cite{Tucker} rigorously justified this geometric model for Lorenz's
original parameters.  Using this model, closed orbits in the Lorenz
dynamical system have been studied combinatorially with symbolic
dynamics on the template.  Indeed, the Lorenz template (see Figure~\ref{F:template}(a))
 can be viewed
as a limit of its periodic orbits, a kind of link with infinitely many
knotted and linked components.  Starting with the template, Birman and Williams \c{BW}
initiated the systematic study of Lorenz knots.  They proved that
infinitely many distinct knot types occur, including all torus knots
and certain cables on torus knots.

Recently, Ghys \c{Gh} established a startling connection with the
periodic orbits in the geodesic flow on the modular surface, which are
in bijection with hyperbolic elements in ${\rm PSL}(2,\Z)$.  Any
hyperbolic matrix $A\in{\rm PSL}(2,\Z)$ defines a periodic orbit,
which Ghys called a {\em modular knot}, in the associated modular flow
on the complement of the trefoil knot in $S^3$.  Ghys proved that
isotopy classes of Lorenz knots and modular knots coincide.  His proof
relies on ingenious deformations that ultimately show that periodic
orbits of the modular flow can be smoothly isotoped onto the Lorenz
template embedded in ${\rm PSL}(2,\R)/{\rm PSL}(2,\Z)$.  (See also
\cite{GhysLeys}, a survey article on this work with breathtaking
images.)

We define {\em Lorenz links} to be all links on the Lorenz template;
i.e., all finite sublinks of the `infinite link' above, a definition
that coincides with Ghys' modular links \c{Gh, Ghysemail}. This
definition is broader than the one used in \c{BW}, which excluded any
link with a parallel cable around any component.  Thus, Lorenz links
are precisely all links as in \cite{BW}, together with any parallel
push-offs on the Lorenz template of any sublinks.  Lorenz knots are
the same in both definitions, but Lorenz links include, for example,
all $(n,n)$--torus links, which are excluded from links in \cite{BW}
for $n\geq 4$. Lorenz braids are all braids on the braid template.

We define T-links as follows.  The link defined by the closure of the
braid $(\s1\cdots\s{r-1})^s$ is a torus link T$(r,s)$.  For $2\leq r_1
\leq \ldots \leq r_k, \ \ 0<s_i, \ \ i=1,\dots,k$, let
T$((r_1,s_1),\ldots,(r_k,s_k))$ be the link defined by the closure of
the following braid, all of whose crossings are positive:
\begin{equation}
\label{E:T--braid}
\TT = (\s1\s2\cdots\s{r_1-1})^{s_1}(\s1\s2\cdots\s{r_2-1})^{s_2}\cdots (\s1\s2\cdots\s{r_k-1})^{s_k}.
\end{equation} 
We call $\TT$ a T--braid, and refer to the link $\T$ that
its closure defines as a T-{\em link}.    See Figure 2 for examples (note that braids are oriented anticlockwise). 
\begin{figure}
\label{F:T-links}
\centerline{\includegraphics[scale=.6]{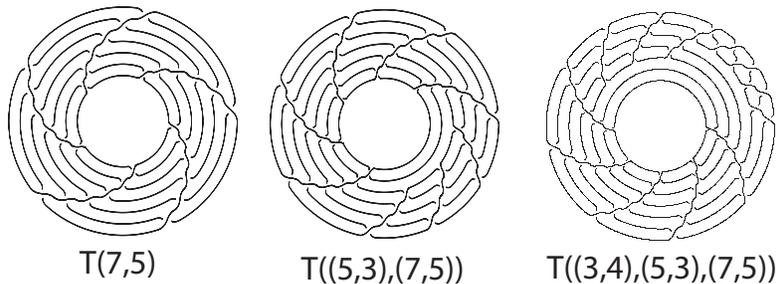}}
\caption{Examples of T-links}
\end{figure}
Note that T--knots, in the case $k=2$, are both more general and less
general than the twisted torus knots studied in \c{CDW, CKP}.  In
those references, twisted torus knots are obtained by performing $s$
{\em full } twists on $r$ strands of a $(p,q)$--torus knot.  This
means that $s_1$ is a multiple of $r_1$, which we do not require in
general.  On the other hand, in those references the twists need not
be positive.

In this paper, in {\bf Theorem~\ref{T:Lorenz-braids and T-braids}} we
establish the following one-to-one correspondence: {\em Every Lorenz
  link is a {\rm T}--link, and every {\rm T}--link is a Lorenz link.}
Among many interesting consequences for both T--links and Lorenz
links, this correspondence suggests a fertile new area for
investigation: the hyperbolic geometry of Lorenz knot complements.

To set up the background, and explain what we learned, recall that
modern knot theory originated with efforts to tabulate knot types.
Starting at the end of the 19th century, tables of knots ordered by
their minimum crossing number were constructed. The early tables went
through 9 and then 10 crossings, and were constructed entirely by
hand.  Roughly 100 years later that project was carried as far as good
sense dictated, when two separate teams, working independently and
making extensive use of modern computers, tabulated all distinct prime
knots of at most 16 crossings, learning in the process that there
there are 1,701,936 of them, now available using the computer program
{\em Knotscape}.  These `knot tables' have served for many years as a
rich set of examples.  The use of minimum crossing number as a measure
of complexity may actually have added to the usefulness of the tables,
because crossing number has limited geometric meaning, so the
tabulated knots serve in some sense as a random collection.

Ghys and Leys \c{GhysLeys} had stressed the scarcity of Lorenz knots
in the knot tables. In particular, Ghys had obtained data which showed
that among the 1,701,936 prime knots with 16 crossings or fewer, only
20 appear as Lorenz knots, with only 7 of those non-torus knots. It
would seem that Lorenz knots are a very strange and unfamiliar
collection.

The study of hyperbolic 3-manifolds, and in particular hyperbolic knot
complements, is a focal point for much recent work in 3-manifold
topology.  Thurston showed that a knot is hyperbolic if it is neither
a torus knot nor a satellite knot. His theorems changed the focus of
knot theory from the properties of diagrams to the geometry of the
complementary space.  Ideal tetrahedra are the natural building blocks
for constructing hyperbolic 3-manifolds, and ideal triangulations can
be studied using the computer program {\em SnapPea}.  There are 6075
noncompact hyperbolic 3-manifolds that can be obtained by gluing the
faces of at most seven ideal tetrahedra \cite{CHW}.  For a hyperbolic knot, the
minimum number of ideal tetrahedra required to construct its
complement is a natural measure of its geometric complexity.

In \c{CDW, CKP}, it was discovered that twisted torus knots occur
frequently in the list of ``simplest hyperbolic knots,'' which are
knots whose complements are in the census of hyperbolic manifolds with
seven or fewer tetrahedra.  Since those twisted torus knots were not
all positive, we collected new data to determine how many were Lorenz
knots.  By the correspondence in Theorem \ref{T:Lorenz-braids and T-braids},
\begin{itemize} 
\item Of the 201 simplest hyperbolic knots, at least 107 are Lorenz knots.  
\end{itemize}
The number 107 could be too small because, among the remaining 94
knots, we were unable to decide whether five of them are or are not
Lorenz.  Many knots in the census had already been identified as
positive twisted torus knots, though their diagrams did not in any way
suggest the Lorenz template.  Lorenz braids for the known Lorenz knots
in the census are provided in a table in Section \ref{S:simplest hyperbolic}.

The data in the census suggests a very interesting question:  
\begin{question}
  Why are so many geometrically simple knots Lorenz knots?  
\end{question}

The heart of the proof of Theorem~\ref{T:Lorenz-braids and T-braids}
is simply that the links in question are all positive, and happen to
have two very different kinds of closed positive braid
representations: as Lorenz braids on the one hand, and as T-braids on
the other hand.
Theorem~\ref{T:Lorenz-braids and T-braids} has immediate consequences
for T-links.  {\bf Corollary~\ref{C:applications1}} asserts that all
of the properties that were established in \cite{BW} for Lorenz links
apply now to T-links, and in particular to positive twisted torus
links. For example, T-links are prime, fibered, non-amphicheiral, and have positive signature.

Another easy consequence for Lorenz links, which was useful for
recognizing Lorenz knots in the census, is {\bf
  Corollary~\ref{C:finite1}}:  Every Lorenz link $\mathcal L$ has
finitely many representations as a Lorenz braid, up to trivial
stabilizations.

Further consequences depend on the observation that the correspondence
in Theorem \ref{T:Lorenz-braids and T-braids} implies certain new
symmetries.  {\bf Corollary~\ref{C:applications2}} applies a somewhat
subtle symmetry of T-links, which generalizes the well-known but not
uninteresting fact that T$(r,s) =$ T$(s,r)$ to Lorenz links.  Going
the other way, there is an obvious symmetry of Lorenz braids by
`turning over the template', which provides a non-obvious involution
for T--braids.  The involution exchanges the total numbers of strands
that are being twisted for numbers of over-passes in the twisted
braid.  See {\bf Corollary~\ref{C:duality}}.  In the special case of
positive twisted torus links, it asserts that T$((r_1,s_1),
(r_2,s_2))$ and T$((s_2, r_2-r_1),(s_1+s_2, r_1))$ have the same link
type.
  
This symmetry, generalized to all T-links below, is quite interesting
in its own regard, and it also enables us to establish new properties
of Lorenz links.  It is a well-known open problem, with many related
important conjectures, to find the precise relationship between the hyperbolic
volume and the Jones polynomial of a knot.  Using
Theorem~\ref{T:Lorenz-braids and T-braids}, the duality of
Corollary~\ref{C:duality}, Thurston's Dehn surgery theorem
\cite{Thurston}, and results in \cite{CK}, we are able to show that
both hyperbolic volume and Mahler measure of Jones polynomials are
bounded for very broadly defined infinite families of Lorenz links.

Let $N>0$.  Let $\mathcal{L}$ be a Lorenz link or T--link
  satisfying any one of the following conditions: The Lorenz braid of
  $\mathcal{L}$ has at most $N$ over-crossing (or under-crossing)
  strands, or equivalently, the T--braid of $\mathcal{L}$ has at most
  $N$ strands (i.e. $r_k\leq N$), or at most $N$ over-passes
  (i.e. $s_1+\ldots+s_k\leq N$).  {\bf Corollaries 5 and 6} assert:
\begin{enumerate}
\item If $\mathcal{L}$ is hyperbolic, its hyperbolic volume is bounded by a constant that depends only on $N$.
\item The Mahler measure of the Jones polynomial of $\mathcal{L}$ is bounded by a constant that depends only on $N$.
\end{enumerate} 
The Jones polynomials of Lorenz links are very atypical, sparse with
small nonzero coefficients, compared with other links of equal
crossing number.  Pierre Dehornoy \cite{Dehornoy} assembled a great
deal of data, but the polynomials were too complicated to pin down
precisely.  {\bf Corollary~\ref{C:Stoimenow}} summarizes the relevant known results about the degrees of
the Jones, HOMFLY and Alexander polynomials of links that can be represented as closed positive braids, and so about Lorenz links.

Continuing our quick review of the paper, we briefly discuss braid
representations of Lorenz links at minimal braid index.  It is known
that the braid index of a Lorenz knot is its {\em trip number}, a
concept that was first encountered in the study of Lorenz knots from
the point of view of symbolic dynamics (see \cite{BW}).
In view of the 1-1 correspondence in Theorem~\ref{T:Lorenz-braids and
  T-braids}, an immediate consequence is that the braid index of each
corresponding T--link is also known.  Nevertheless, a problem arises: If
$\mathcal T$ is a T--link, the trip number of its Lorenz companion
$\mathcal L$ is not easily computed from the defining parameters for
the T--link.  In {\bf Corollary~\ref{C:braid index of T--links}}, we
give an explicit formula for computing it directly from the sequence of
integer pairs $((r_1,s_1),\dots,(r_k,s_k))$ that define $\mathcal T$.

In $\S$\ref{S:Minimal braid index representations}, we prove {\bf
  Theorem~\ref{T:Lorenz-braids and t-braids}}, which establishes for
any Lorenz link $\mathcal{L}$, a correspondence between its Lorenz
braid representations and particular factorizations of braid words in
the braid group $B_t$, where $t$ is the minimal braid index of
$\mathcal{L}$.  This theorem is a strong form of Proposition 5.1 of
\cite{BW}, and is interesting because it applies to T--links as well
as to Lorenz links.  

We return to the hyperbolicity question: When
is a Lorenz link hyperbolic?  We could not answer that question, but
as a starter {\bf Corollary~\ref{C:algorithm}} gives a fast algorithm
to decide when a Lorenz link is a torus link.

Here is a guide to this paper.  In $\S$\ref{S:preliminaries} we set up
our notation and prove a basic lemma about the repeated removal of
trivial loops in a Lorenz braid.  The lemma will be used in the proofs
of Theorems~\ref{T:Lorenz-braids and T-braids} and
\ref{T:Lorenz-braids and t-braids}.  In $\S$\ref{S:Lorenz links and
  T-links}, we state and prove Theorem~\ref{T:Lorenz-braids and
  T-braids} and Corollaries~\ref{C:applications1}--\ref{C:braid index
  of T--links}.  In $\S$\ref{S:Minimal braid index representations} we
prove Theorem~\ref{T:Lorenz-braids and t-braids}
and Corollary~\ref{C:algorithm}.  In $\S$\ref{S:simplest
  hyperbolic}, we discuss and provide Lorenz data for the simplest
hyperbolic knots.  Open questions are scattered throughout the paper.

{\bf Acknowledgments} We thank the two referees of an early version of
this paper for their thoughtful and constructive comments.  We thank
Slavik Jablan for sharing with us the list of Lorenz knots up to 49
crossings, and verifying for us that certain simplest hyperbolic knots
are not on this list, as discussed in Section \ref{S:simplest
  hyperbolic}.  We thank Volker Gebhardt for sharing his program to
draw closed braids, which was used to make Figure~\ref{F:T-links}.  We
thank Pierre Dehornoy, Etienne Ghys, Juan Gonzalez-Meneses, Slavik
Jablan, Michael Sullivan, and Robert Williams for carefully answering
our questions.

\section{Preliminaries}
\label{S:preliminaries}  

We defined a Lorenz link to be any finite collection of closed orbits
on the Lorenz template, which supports a semiflow.  The template is a
branched 2-manifold embedded in $\mathbb{R}^3$, as illustrated in
Figure~\ref{F:template}.
In the right sketch, the Lorenz template has been cut open to give a
related template for Lorenz braids, which inherit an orientation from
the template, top to bottom.  The crossings in Figure~\ref{F:lorenz
  braid} will be called {\em positive} crossings.  Although this
convention is opposite to the usual one in knot theory, it matches
\cite{BW} and has often appeared in the related literature, so we
continue to use it now.

\begin{example} \label{Ex:favorite} \rm
Figure~\ref{F:lorenz braid} gives an example of a Lorenz braid. It
becomes a Lorenz knot after connecting the strands as in a closed
braid, on the template.  This example will be used throughout the
paper to illustrate our ideas, so Figure~\ref{F:lorenz braid}
contains features that will be explained later.
\begin{figure}[htpb!]
\centerline{\includegraphics[scale=.6] {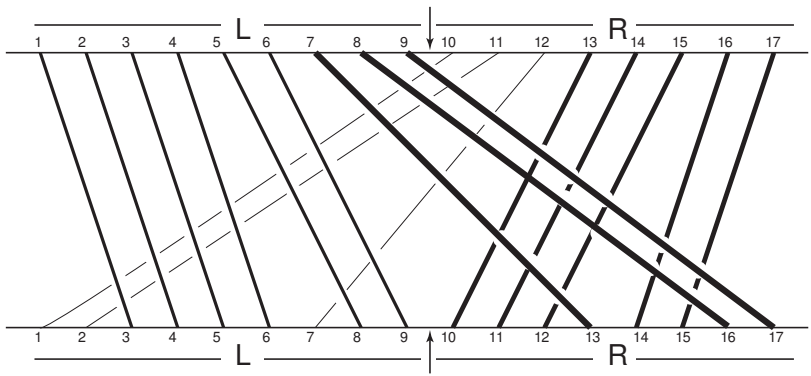}}
\caption{$\vec{d}_L = \langle 2, 2, 2, 2, 3, 3, 6 ,8, 8 \rangle = \langle 2^4,3^2,6^1,8^2\rangle$}
\label{F:lorenz braid}
\end{figure}
\end{example}

The Lorenz braid $\mathbb{L}$ is determined entirely by its
permutation, because any two strands cross at most once.  In a Lorenz
braid two overcrossing (resp. undercrossing) strands never intersect,
so the permutation associated to the overcrossing strands uniquely
determines the rest of the permutation.  Therefore, $\mathbb{L}$ is
determined by just the permutation associated to its overcrossing
strands.

Assume there are $p>1$ overcrossing strands.  On each overcrossing
strand the position of the endpoint will always be bigger than that of
the initial point.  Suppose the $i^{\rm th}$ strand begins at $i$ and
ends at $i+d_i$.  Since two overcrossing strands never cross, we have
the following sequence of positive integers: 
$$ d_1 \leq d_2\leq \cdots \leq d_{p-1}\leq d_p $$ 

Lorenz braids that have unknotted closure were classified in Corollary
5.3 of \cite{BW}.  Excluding the two trivial loops that are parallel
to the two boundary components, it was proved that a Lorenz knot is
unknotted if and only if the following condition holds:
$$ i + d_i > p \quad {\rm if\ and\ only\ if} \quad i=p $$ 
Since $1 \leq d_1 \leq d_2 \leq \cdots \leq d_p$, it follows that $d_i
= 1$ for every $i\leq p-1$.  But $d_p$ can be $1,2,3,\dots$, and
these are the only ways to obtain the unknot.

In view of this classification, we can make two assumptions: (i)
$2\leq d_1$ and (ii) $d_{p-1} = d_p$.  Otherwise, if $d_1=1$ then
$\L=\sigma_1\L'$ for a Lorenz braid $\L'$ on the last $n-1$ strands,
so that $\L$ can be trivially destabilized on its left side.  We get
similarly trivial destabilizations on the right if $d_{p-1} < d_p$.
As we have seen, the only closed orbits omitted by making these
assumptions are the Lorenz unknots.

We collect this data in the following vector (see \cite{BW}):
\begin{equation} \label{E:1}
\vec{d}_L  = \langle d_1,\dots,d_p\rangle, \quad  2\leq d_1, \ d_{p-1} = d_p, \ {\rm and\ each}\ d_i\leq d_{i+1}
\end{equation}
The vector $\vec{d}_L$ determines the positions of the L 
(overcrossing) strands. The R (undercrossing) strands fill in the
remaining positions, in such a way that all crossings are L-strands
crossing over R-strands.  In Figure \ref{F:lorenz braid}, the arrows
separate the left and right strands.  Each $d_i$ with $i=1,\dots,p$ is
the difference between the initial and final positions of the
$i^{\rm{th}}$ overcrossing strand.  The integer $d_i$ is also the
number of strands that pass under the $i^{\rm{th}}$ braid strand.  The
vector $\vec{d}_L$ determines a closed braid $\L$ on $n = (p+d_p)$
strands, which we call the {\em Lorenz braid} representation of the
Lorenz link $\cal L$.  All non-trivial Lorenz links arise in this way.

The overcrossing strands travel in groups of {\em parallel} strands,
which are strands of the same slope, or equivalently strands whose
associated $d_i$'s coincide.  If $d_{\mu_j} = d_{\mu_j + 1} = \cdots =
d_{\mu_j + s_j -1}$, where $s_j$ is the number of strands in the
$j^{\rm{th}}$ group, then let $r_j = d_{\mu_j}$.  Thus, we can write
$\vec{d}_L$ in the form:
\begin{equation}
\label{E:2}
\vec{d}_L =  \langle d_{\mu_1}^{s_1},\dots, d_{\mu_k}^{s_k} \rangle =  \langle r_1^{s_1}, \dots, r_k^{s_k}\rangle, \ \  1\leq s_i, \ \ 2\leq r_1,s_k, \ \ {\rm and} \ \  r_i < r_{i+1}.
\end{equation} 
where $r_i^{s_i}$ means $r_i,\ldots, r_i$ repeated $s_i$ times.  
Note that
$$ p=s_1+\cdots + s_k, \quad d_1 = r_1, \quad d_p = r_k $$
The {\em trip number} $t$ of a Lorenz link $\mathcal L$ is given by
\begin{equation}
\label{E:trip number}
t = \#\{ \ i \ | \  i + d_i > p \ {\rm where}  \ 1\leq i\leq p \ \}
\end{equation}
The trip number is the minimum braid index of $\mathcal L$, a fact which was conjectured in \cite{BW} and proved in \cite{FW}.  

In Example~\ref{Ex:favorite}, $\langle 2, 2, 2, 2, 3, 3, 6, 8, 8
\rangle = \langle 2^4,3^2,6^1,8^2\rangle$.  Thus, $p=9,\ k=4,\ r_k=8$,
and $n= p + r_k = 17$ is the braid index of $\mathbb{L}$.  The trip
number $t=3$ is the braid index of the Lorenz link given by the
closure of $\mathbb{L}$.

Our first new result is little more than a careful examination of the
proof of Theorem 5.1 of \cite{BW}.  This lemma will be used in the
proofs of Theorems~\ref{T:Lorenz-braids and T-braids} and
\ref{T:Lorenz-braids and t-braids} of this paper.
\begin{lemma} 
\label{L:Markov sequence}
Let $\mathbb{L}$ be a Lorenz braid defined by $\vec{d}_L = \langle
r_1^{s_1},\dots,r_k^{s_k} \rangle = \langle d_1,\dots,d_p\rangle$, so
$\mathbb{L}$ is a braid on $p+r_k$ strands.  Then there is a sequence of
closed positive braids:
$$ \mathbb{L} = \mathbb{L}_0 \to \mathbb{L}_1\to \mathbb{L}_2\to \cdots \to \mathbb{L}_p \to \cdots \to\mathbb L_{p+r_k-t}$$ 
where each $\mathbb{L}_{i+1}$ is obtained from $\mathbb{L}_i$ by
single move that reduces the braid index and also the crossing
number by $1$.  Each $\mathbb{L}_i$ represents the same Lorenz link $\mathcal{L}$.
The intermediate braid $\mathbb{L}_p$ in the sequence has $r_k$ strands, and the final braid $\mathbb L_{p+r_k-t}$ has $t$ strands, which is the minimum braid index of $\mathcal L$.   
\end{lemma}

\pf  We will find the required Markov sequence by using a geometric trick that was introduced in \cite{BW}.    

The letters $1,2,\dots,p$ in a Lorenz permutation are said to be in
the {\em left} or L-group, and the letters $p+1,\dots, p+ r_k$ are in
the {\em right} or R-group. Each strand in a Lorenz braid begins and
ends at a point which is either in L or in R, therefore the strands
divide naturally into 4 groups: strands of type LL, LR, RL and RR,
where strands of type LL (resp. LR) begin on L and end on L (resp. R),
and similarly for types RL and RR.  In the example in Figure
\ref{F:lorenz braid}, those of type LR are the thickest, with those of
type RR, LL and RL each a little thinner than their
predecessors.

\begin{figure}[htpb!]
\centerline{\includegraphics[scale=.55] {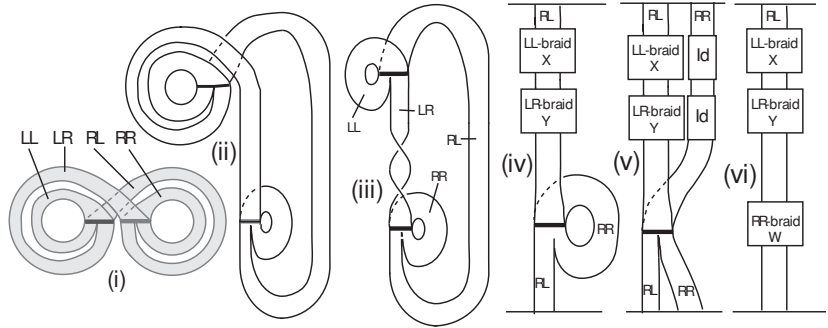}}
\caption{Cutting open the Lorenz template}
\label{F:cutting-open}
\end{figure}

By definition of the trip number $t$ in (\ref{E:trip number}), there are  $p-t, \ t, \ t$ and $r_k - t$ strands of type LL,LR,RL and RR respectively. 
In sketch (i) of Figure~\ref{F:cutting-open} we have cut open the Lorenz template, snipping it open between two orbits, as was done in \cite{BW}, so that the template itself divides naturally into bands of type LL,LR,RL and RR. 

In sketch (ii), we stretched out the band that contains all of the
strands of type LR, and in sketch (iii) we uncoiled that band,
introducing a full twist into the $t$ strands of type LR.  This
uncoiling can be regarded as having been done one strand at a time,
and when we think of it that way, it becomes a sequence of $t$ moves,
each reducing the braid index by 1.  Observe that when we `uncoil' the
outermost arc (and the ones that follow too), we trade one arc in the
braid $\mathbb{L}_i$ for a `shorter' arc in the braid
$\mathbb{L}_{i+1}$, reducing the braid index by 1.  This process has
been repeated $t$ times in the passage from sketch (ii) to sketch
(iii), because there are $t$ strands in the LR braid.  The uncoiling
takes positive braids to positive braids, although the property of
being a {\em Lorenz braid} is not preserved.  After the $t$ Markov
moves illustrated in the passage to sketch (iii), the braid index will
have been reduced from $p+r_k$ to $p+r_k -t$.

\begin{figure}[htpb!]
\centerline{\includegraphics[scale=.6] {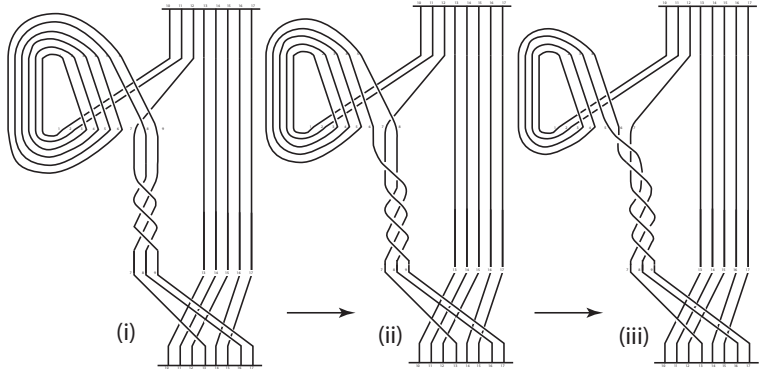}}
\caption{Uncoiling the LL braid: an example}
\label{F:uncoiling-example}
\end{figure}

We turn our attention to the LL subbraid in sketch (iii), which has $p
- t$ strands.  The strands of type RL and of type LR both have $t$
strands. From this it follows that when the LL band is uncoiled, we
obtain a subbraid on $t$ strands that joins the RL subbraid to the LR
subbraid, as illustrated in sketch (iv) of
Figure~\ref{F:cutting-open}.  Figure~\ref{F:uncoiling-example},
illustrates via an example the uncoiling that leads to the LL braid.
The example shown in Figure~\ref{F:uncoiling-example} is the
Lorenz braid from Example \ref{Ex:favorite}, shown in Figure~\ref{F:lorenz braid}.
It's a rather simple example because the trip number $t=3$.  Sketch
(i) in Figure~\ref{F:uncoiling-example} corresponds to sketch (iii) in
Figure~\ref{F:cutting-open}. Sketches (ii) and (iii) of
Figure~\ref{F:uncoiling-example} show two destabilizations, and
correspond to two steps in the passage from sketch (iii) to sketch
(iv) of Figure~\ref{F:cutting-open}.  

In the LL braid, we uncoil the $i^{th}$ strand, which is the outer
coiled strand in the leftmost sketch in
Figure~\ref{F:uncoiling-example}.  Let $\gamma_i$ be this outer arc
running clockwise and crossing over three strands in the sketch.  If
$\pi$ is the permutation associated to the Lorenz braid $\mathbb{L}$,
say $\gamma_i(0)$ corresponds to the bottom endpoint $i$,
$\gamma_i(1/2)$ corresponds to the top endpoint $i$, and $\gamma_i(1)$
corresponds to $\pi(i)=i+d_i$.  Therefore, in the middle sketch of
Figure \ref{F:uncoiling-example}, $\gamma_i$ is replaced by an arc
that contributes $(\s 1 \s 2 \cdots \s {d_i-1}) $ to the LL braid.
Continuing in this way, the LL braid determines $\prod_{i=1}^{i=p-t} (\s 1 \s 2 \cdots \s {d_i-1})$.

We return to Figure~\ref{F:cutting-open}.  There are $p-t$ strands in
the LL braid, so there are $p-t$ strands that are uncoiled in sketch
(iv). The braid index will go down from $(p+r_k-t)$, in sketch (iii),
to $(p+r_k-t) - (p-t) = r_k$ in sketch (iv).  In sketch (v) we have
cut open the right coil, exhibiting the braid template for the
$r_k$--braid to be studied in Theorem~\ref{T:Lorenz-braids and
  T-braids}.  This braid was not considered in \cite{BW}.  

The final braid in the sequence is illustrated in sketch (vi). It is
obtained from the braid in sketch (iv) by uncoiling the RR braid.  It
was proved in \cite{BW} to have braid index $t$, where $t$ is the trip
number.  We will study it further in $\S$\ref{S:Minimal braid index
  representations}.

Each braid $\mathbb{L}_{i+1}$ in the sequence that we just described has braid index one less than that of its predecessor $\mathbb{L}_i$.  To prove the assertion about crossing number, observe that since each $\mathbb{L}_i$ is a positive braid, an Euler characteristic count shows that  
$$2g(\mathcal{L}) = c(\mathbb{L}_i)-n(\mathbb{L}_i) - \mu(\mathcal{L}) + 2,$$ 
where $c$ is the crossing number of the positive braid $\mathbb{L}_i$, $n$ is its braid index, $g$ is the genus, and $\mu$ is the number
of components of $\mathcal{L}$.  But then $c(\mathbb{L}_i)-n(\mathbb{L}_i)$ is a topological invariant of $\mathcal{L}$, so when $n$ is reduced $c$ must be reduced too.
This completes the proof of Lemma \ref{L:Markov sequence}.
\done  

\section{Lorenz links and T-links}
\label{S:Lorenz links and T-links}

T-links  were defined  in $\S$\ref{S:introduction} above. 
\begin{theorem}\label{T:Lorenz-braids and T-braids}
Every Lorenz link is a {\rm T}--link, and every {\rm T}--link is a Lorenz link.  

Precisely, if a link $\mathcal{L}$ is represented by a Lorenz braid
$\mathbb{L}$ on $p+r_k$ strands with $\vec{d}_L = \langle r_1^{s_1},
\dots, r_k^{s_k}\rangle$, then $\mathcal{L}$ also has an $r_k$--braid
representation $\mathbb{T}$, given in {\rm (\ref{E:T--braid})}, which
exhibits it as a {\rm T}--link. Moreover, every {\rm T}--link arises
in this way from some Lorenz link.
\end{theorem}
\pf
We begin by introducing convenient notation.
Let $v,w$ be positive integers with $v<w$. 
Let $[v,w]=\sigma_v\sigma_{v+1}\cdots\sigma_{w-1}$ and let $[w,v]=\sigma_{w-1}\cdots\sigma_{v+1}\sigma_v$.
If $u<v<w$, we have a very simple {\it product rule}:
\begin{equation}
\label{E:product rule}
[u,v][v,w] = [u,w]
\end{equation}
In the braid group $B_n$, an {\em index shift relation} holds:
$$ (\s1\s2\dots\s {k-1}\s k)(\s j) = ( \s {j+1})(\s 1\s 2\dots\s {k-1}\s k), \ \ j=1,2,\dots,k-1 $$
The index shift relation can be expressed in our new notation as:
\begin{equation}
\label{E:index shift using brackets}
[1,w][u,v] = [u+1,v+1] [1,w]\ \ {\rm if} \ v<w
\end{equation}

We now study the braids in sketch (v) of Figure~\ref{F:cutting-open} in detail.  
In the proof of Lemma~\ref{L:Markov sequence}, which was based upon
Figure~\ref{F:cutting-open}, we produced a Markov sequence from our
original braid $\mathbb{L}\in B_{p+r_k}$ to a braid $\mathbb{L}_p \in
B_{r_k}$.  It is clear from sketch (v) of
Figure~\ref{F:cutting-open} that $\mathbb{L}_p$ is a product of braids
$\mathbb{XYZ}\in B_{r_k}$, where $\mathbb{X}$ comes from the LL braid,
$\mathbb{Y}$ from the LR braid, and $\mathbb{Z}$ from the RL and RR braids.  Both $\mathbb{X}$ and $\mathbb{Y}$ use
only the first $t$ strands, their remaining $r_k-t$ strands being the
identity braid, but $\mathbb{Z}$ uses all $r_k$ strands.

In the proof of Lemma 1, we learned that the braid word that describes $\mathbb{X}$ is a product of the form: 
\begin{equation}
\label{E:expression for X}
\mathbb{X} = \prod_{i=1}^{p-t} (\s 1 \s 2 \cdots \s {d_i-1}) = \prod_{i=1}^{p-t}\; [1,d_i]
\end{equation}

We will identify the braid $\mathbb{YZ}$ that is associated to the LR,
RR and RL subbraids in Figure~\ref{F:cutting-open}.  From sketch (iii)
of Figure~\ref{F:cutting-open}, one sees immediately that $\mathbb{Y}
= [1,t]^t$.  Let $\mathbb{Z}_t$ be the braid on $r_k$ strands that is
created when the strand that begins at position $t$ and ends in
position $(t+r_k) -t = r_k$ crosses over all the intermediate strands,
with every strand that is not crossed remaining fixed.  So
$\mathbb{Z}_t=[t,r_k]$.  Let $\mathbb{Z}_{t-i}$ be the braid that is
associated to the strand that begins at $t-i$, where
$i=0,1,\dots,t-1.$ This strand crosses over all the intermediate
strands, but all strands that are not crossed remain fixed.  Therefore
$$\mathbb{Z}_{t-i} = [t-i, t-i + d_{p-i} - t] = [t-i,d_{p-i}-i].$$ 
Our next claim is the key to the proof of Theorem~\ref{T:Lorenz-braids and T-braids}:

Let $\mathbb{Y} = [1,t]^t$ and let $\mathbb{Z} = \mathbb Z_t \mathbb Z_{t-1}\cdots \mathbb Z_1$. We claim:
\begin{eqnarray}
\label{E:YZ1}
\mathbb{Y}\mathbb{Z}  =  [1,d_{p-t+1}][1,d_{p-t+2}]\cdots [1, d_{p}]
\end{eqnarray}
We will prove (\ref{E:YZ1}) by induction on $i$, where $i=0,1,\dots,t-1$.
Let $\mathbb Y_i=[1,t]^{i+1}$, and let $\mathbb Z(i) = \mathbb
Z_t\cdots \mathbb Z_{t-i}$.
If $i=0$, we have: $\mathbb Y_0\mathbb{Z}(0) = \mathbb Y_0\mathbb
Z_t = [1,t][t,d_p] = [1,d_p]$, so the induction begins.  Choose any
$i$ with $0<i<t$ and assume, inductively, that
$$\mathbb{Y}_i\mathbb{Z}_t\mathbb{ Z}_{t-1}\dots \mathbb{Z}_{t-i} = [1,d_{p-i}]\cdots [1, d_{p-1}][1,d_p].$$ 
Since $\mathbb{Z}_{t-(i+1)} = [t-(i+1),d_{p-(i+1)}-(i+1)]$, by our induction hypothesis,
\begin{equation}
\label{E:to be simplified}
\mathbb{Y}_{i+1}\mathbb{Z}_t\dots\mathbb{Z}_{t-(i+1)}= [1,t][1,d_{p-i}]\cdots[1,d_p][t-(i+1),d_{p-(i+1)}-(i+1)]
\end{equation}
We must prove that the right-hand side of (\ref{E:to be simplified})
equals $[1, d_{p-(i+1)})]\cdots[1,d_p]$.  This exercise
reveals some subtle consequences of the braid relations.  The
product rule (\ref{E:product rule}) and index shift relation
(\ref{E:index shift using brackets}) will play crucial roles.  We
claim that, as a consequence of (\ref{E:index shift using brackets}),
the factor $[t-(i+1),d_{p-i}-(i+1)]$ on the right in (\ref{E:to be
  simplified}) can be shifted to the left, past all but one of the
brackets on its left, changing its name as it does so. The reasons
are:
\begin{itemize}
\item By our basic definition of the $d_i$'s, we know $d_{p-i}\leq d_{p-i-j}$ for all $j>0$. 
\item Every strand of type LR crosses over all strands of type RL. Since there are $t$ strands of type RL, we conclude that $t<d_{p-t+1}\leq d_{p-i}$ for all $i=0,1,\dots,t-1$. 
\end{itemize}
 These imply that (\ref{E:index shift using brackets}) is applicable $i+1$ times, so that (\ref{E:to be simplified}) simplifies as follows:  
\begin{eqnarray}
\label{E:braid relations}
\nonumber
\mathbb{Y}_{i+1}\mathbb{Z}_t\dots\mathbb{Z}_{t-(i+1)} &=&[1,t][1,d_{p-i}]\cdots [1, d_{p-1}][1,d_p][t-(i+1),d_{p-i}-(i+1)] \\
\nonumber
&=&[1,t][1,d_{p-i}]\cdots [1, d_{p-1}][t-i,d_{p-i}-i][1,d_p] = \cdots = \\
&=&[1,t][t,d_{p-(i+1)})][1,d_{p-i}]\cdots [1, d_{p-1}][1,d_p]
\end{eqnarray}
Finally, we use (\ref{E:product rule}) to combine the two leftmost terms in (\ref{E:braid relations}),
obtaining:
$$ = [1, d_{p-(i+1)}][1,d_{p-i}]\cdots [1, d_{p-1}][1,d_p]$$
This is the desired expression, so (\ref{E:YZ1}) is proved.

Let's put together the expression for
$\mathbb{X}$ in (\ref{E:expression for X}) and $\mathbb{YZ}$ in
(\ref{E:YZ1}).  After collecting like terms in the previous
expression, we obtain:
\begin{eqnarray}
\label{E:XYZ1}
\nonumber
\mathbb{XYZ}  & = & [1,d_1] [1,d_2] \cdots [1,d_{p-t}] [1,d_{p-t+1}] \cdots [1,d_p] \\
&=& [1,d_1][1,d_2]\cdots [1,d_p] \\
\label{E:XYZ2}
&=& [1,r_1]^{s_1}\cdots  [1, r_k]^{s_k}  =  (\s1\cdots\s{r_1-1})^{s_1}\cdots (\s1\cdots\s{r_k-1})^{s_k} 
\end{eqnarray}
where, in the passage (\ref{E:XYZ1}) $\to$ (\ref{E:XYZ2}), we have collected those terms for which successive entries $d_i$ and $d_{i+1}$ coincide, as in the passage (\ref{E:1}) $\to$ (\ref{E:2}).  
But (\ref{E:XYZ2}) is precisely what we claim in Theorem~\ref{T:Lorenz-braids and T-braids}. 

The only remaining question is whether every T--braid 
is obtained from some Lorenz link.  Suppose we are given an arbitrary T--braid $\TT$, whose closure is 
the T--link T$((r_1,s_1),\ldots,(r_k,s_k))$.  Let $\vec{d}_L = \langle r_1^{s_1}, \dots, r_k^{s_k}\rangle$, 
which determines a Lorenz braid $\L$.
By the proof above, $\L$ is braid-equivalent to $\TT$.
This completes the proof of Theorem~\ref{T:Lorenz-braids and T-braids}.  
\done  

There are many consequences of Theorem~\ref{T:Lorenz-braids and T-braids}. 
We can immediately establish many new properties for T--links:
\begin{corollary}
\label{C:applications1} 
The following properties of Lorenz links are also satisfied by
all T-links, and so in particular by positive twisted torus links.
\begin{itemize}
\item[{\rm (i)}] T--links are prime. 

\item[{\rm (ii)}] T--links are fibered. Their genus $g$ is given by the formula $2g = c - n + 2 -\mu$, where $n$ is the braid index of any positive braid representation, $c$ is the crossing number of same, and $\mu$ is the number of components.   

\item[{\rm (iii)}] T--links are non-amphicheiral and have positive signature.
\end{itemize}
\end{corollary}

\pf Property (i) was proved by Williams in \cite{W}.  His proof is
interesting for us, because it illustrates the non-triviality of
Theorem~\ref{T:Lorenz-braids and T-braids}.  Williams used the fact
that all Lorenz links embed in the Lorenz template, and if a Lorenz
link was not prime then a splitting 2-sphere would have to intersect
the template in a way that he shows is impossible.  Without the
structure provided by the template, it seems very difficult to
establish this result for T--links.  Properties (ii) and (iii) were
established in \cite{BW} for Lorenz links.  By
Theorem~\ref{T:Lorenz-braids and T-braids}, they also hold for
T-links.  \done

Theorem~\ref{T:Lorenz-braids and T-braids} provides an easy proof that there are
finitely many non-trivial Lorenz braid representations for any Lorenz
link:

\begin{corollary}\label{C:finite1} 
Every Lorenz link $\mathcal{L}$ has finitely many Lorenz braid representations up to trivial stabilizations.
\end{corollary}
\pf 
By Theorem~\ref{T:Lorenz-braids and T-braids}, there is a one-to-one correspondence
between Lorenz braid representations and T--braid representations of
$\mathcal{L}$.  For any T--braid representation, $c,\,g,\,n,\,\mu$ as
above satisfy $2g = c - n + 2 -\mu$.  From (\ref{E:2}), $s_k \geq 2$, which implies that $c\geq
2(n-1)$.  Hence, $n \leq c-n+2 = 2g+\mu$.  (Since $n=2g+\mu$ for
T$(2,n)$, this inequality is sharp.)  Therefore,
\begin{equation}\label{E:c4g}
c = 2g +\mu + n -2 \leq 4g + 2\mu - 2
\end{equation}
With $c$ bounded, there are only finitely many T--braid representations of
$\mathcal{L}$.  \done

Although any given Lorenz knot appears infinitely often as a component
in its many parallel copies, this does not contradict Corollary
\ref{C:finite1}.  Essentially, we are counting links
rather than their individual components.
Parallel push-offs result in distinct Lorenz links because the Lorenz
template has non-trivial framing, so any two such parallel components
are non-trivially linked.  For example, parallel push-offs of the
unknot are $(n,n)$--torus links.

\subsection{Symmetries}  Another application of the correspondence in
Theorem~\ref{T:Lorenz-braids and T-braids} is to exploit natural
symmetries on one side to establish unexpected equivalences on the
other.  We will do this in both directions; first, from T--links to
Lorenz links:

\begin{corollary}\label{C:applications2} 
Let $\mathbb{L}_1$ be the Lorenz braid defined by $\vec{d}_{L_1} = \langle
r_1^{s_1},\ldots, r_k^{s_k}\rangle$ such that $r_{k-1}\leq s_k$ and
$s_i=n_i\cdot r_i$, for any positive integers $n_i$, with
$i=1,\ldots,k-1$.  Let $\mathbb{L}_2$ be the Lorenz braid defined by
$\vec{d}_{L_2} = \langle r_1^{s_1},\ldots,r_{k-1}^{s_{k-1}},\, s_k^{r_k}\rangle$.  
Then $\mathbb{L}_1$ and $\mathbb{L}_2$ both represent the same Lorenz link.
\end{corollary}
\pf 
By Theorem \ref{T:Lorenz-braids and T-braids}, the closure of $\mathbb{L}_1$ is
$\mathcal{T}_1=$ T$((r_1,s_1),\ldots,(r_k,s_k))$, and the closure of
$\mathbb{L}_2$ is $\mathcal{T}_2=$ T$((r_1,s_1),\ldots,(r_{k-1},s_{k-1}), (s_k,r_k))$.  We claim
that $\mathcal{T}_1$ and $\mathcal{T}_2$ are isotopic, so $\mathbb{L}_1$
and $\mathbb{L}_2$ both represent the same link.

\begin{figure}[htpb!]
\centerline{\includegraphics[scale=.4] {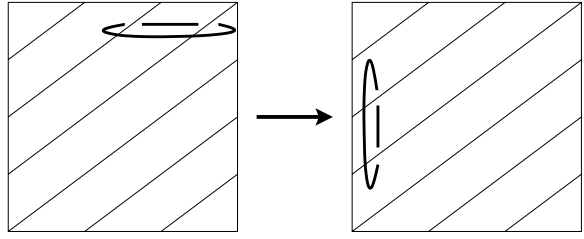}}
\caption{Isotopy in the proof of Corollary \ref{C:applications2}}
\label{F:torus}
\end{figure}

For each $i=1,\ldots,k-1$, the isotopy is the same as in the proof of Lemma 3.1.1 of \cite{De2}.
For all $i,\ s_i=n_i\cdot r_i$, so $\mathcal{T}_1$ is obtained by $1/n_i$--Dehn surgeries
on a nested sequence of unknots $U_i$ that encircle $r_i$ strands of an $(r_k,s_k)$--torus
link.  To obtain $\mathcal{T}_2$, we isotope the $(r_k,s_k)$--torus
link to a $(s_k,r_k)$--torus link.  Since for all $i$, $r_i<r_k$ and $r_i\leq s_k$, we
can slide all $U_i$ along the torus link from the meridinal to the
longitudinal direction (see Figure \ref{F:torus}), and then perform
the same Dehn surgeries.  
\done

\begin{example} \label{E:torus} 
For $k=2$, Corollary \ref{C:applications2} implies that a non-trivial (i.e. $2\leq r_1 < r_2$ and
$1\leq s_1,\ 2\leq s_2$) twisted torus link is in fact a torus link if $r_1=s_2$ and $s_1=n\cdot r_1$.
For example, the Lorenz braid given by $\vec{d}_L = \langle 3^6, 8^3
\rangle$ represents a torus knot:
$$ T((3,6),(8,3)) = T((3,6),(3,8)) = T(3,14) $$
In addition, because the Lorenz braids given by $\vec{d}_L = \langle 3^6, 8^3
\rangle$ and $\langle 3^{14} \rangle$ represent the same link, we see
that the integer $k$ is not an invariant of link type.
\end{example}

Going in the opposite direction, a natural symmetry of Lorenz links provides a
far-reaching application of Theorem \ref{T:Lorenz-braids and T-braids}
to T-links.  Observe that a rotation of $\pi$ about the $z$ axis in
Figure~\ref{F:template} is a symmetry of the Lorenz template. If a
Lorenz braid $\overline{\L}$ is obtained from a Lorenz braid $\L$ by
this rotation, we will say that $\L$ and $\overline{\L}$ are {\em dual
  Lorenz braids}.  Both $\L$ and $\overline{\L}$ determine the same
Lorenz link $\mathcal{L}$.  By Theorem~\ref{T:Lorenz-braids and
  T-braids}, there is a corresponding duality between T--braid
representations of $\mathcal{L}$, which vastly generalizes the
well-known duality for torus links, T$(r,s)=$T$(s,r)$.
\begin{corollary}
\label{C:duality}\quad   Let
\begin{equation}
\label{E:duality 1}
\overline{r}_1 = s_k, \ \ \overline{r}_2 = s_k+s_{k-1},\ \ \dots, \ \ \overline{r}_k = s_k+s_{k-1}+\cdots + s_1
\end{equation}
\begin{equation}
\label{E:duality 2}
\overline{s}_1 = r_k - r_{k-1}, \ \ \overline{s}_2 = r_{k-1}-r_{k-2},\ \ \dots,\ \ \overline{s}_{k-1} = r_2 - r_1, \ \ \overline{s}_k = r_1.
\end{equation}  Then
  T$((r_1,s_1),\ldots,(r_k,s_k))$ \ and \ 
  T$((\overline{r}_1,\overline{s}_1),\ldots,(\overline{r}_k,\overline{s}_k))$
  have the same link type.
\end{corollary}
\pf
By Theorem~\ref{T:Lorenz-braids and T-braids}, the claim can be proved for a pair of dual Lorenz braids.
Let $\mathcal{L}$ be a Lorenz link with a Lorenz braid representation $\L$, defined by
$$\vec{d}_L = \langle d_1, d_2, \dots, d_p \rangle =
\langle r_1^{s_1},r_2^{s_2},\dots,r_k^{s_k}\rangle $$ 
We claim that 
$\mathcal{L}$ also has a dual Lorenz braid representation
$\overline{\mathbb{L}}$, defined by the {\em dual vector}
$$\overline{\vec{d}_L} = \langle \overline{d}_1, \overline{d}_2, \dots, \overline{d}_{\overline{p}} \rangle =
\langle\overline{r}_1^{\overline{s}_1},\overline{r}_2^{\overline{s}_2},\dots,\overline{r}_k^{\overline{s}_k}\rangle$$
with $\overline{p} = d_p$ and $\overline{d}_{\overline{p}} = p$, and
the rest of the entries determined by $\vec{d}_L$ as in (\ref{E:duality 1}) and (\ref{E:duality 2}). 
This would imply that $\L$ and $\overline{\L}$ have the same braid
index, $p+d_p = \overline{p} + \overline{d}_{\overline{p}}$, and both
represent $\mathcal{L}$.  This correspondence is equivalent to the
statement in the theorem: for every T--braid representation of
$\mathcal{L}$ in $B_{r_k}$, there is a dual T--braid representation of
$\mathcal{L}$ in $B_p$, given by (\ref{E:duality 1}) and (\ref{E:duality 2}).

The strands in $\L$ divide into overcrossing strands and undercrossing strands, whose roles are interchanged when we pass from $\L$ to $\overline{\L}$.  See Figure~\ref{F:lorenz braid}.  From this it follows that $\overline{p} = d_p$ and $\overline{d}_{\overline{p}} = p.$ But
then $p+d_p = \overline{p} + \overline{d}_{\overline{p}}$, so that both have the same braid
index.  The dual braid is then simply the
original one, flipped over so that strand $i$ becomes strand
$p+d_p-i+1$. Clearly both determine the same link $\mathcal{L}$.

A {\em crossing point} in the braid $\L$ or $\overline{\L}$ means a double point in the projected image.  Two overcrossing strands in $\L$, (and also in $\overline{\L}$) are said to be {\em parallel} when they contain the same number of crossing points.  Observe that the overcrossing strands in $\L$ divide naturally into packets of parallel strands, where the $i^{\rm th}$ group of parallel strands contains $s_i$ strands, each of which has $r_i$ crossing points.  In the same way, there is a different subdivision of the overcrossing strands of $\overline{\L}$, with the  $j^{\rm th}$ group of parallel strands containing $\overline{s}_j$ strands, each having $\overline{r}_j$ crossings. 

Now observe that there are blank spaces between the endpoints of the $i^{\rm th}$ and $(i+1)^{\rm st}$ group of overcrossing strands in $\L$ for exactly $r_{i+1} - r_i$ overcrossing strands of $\overline{\L}$.     Taking into account that strand $i$ in $\L$ becomes strand $n-i+1$ in  $\overline{\L}$, it follows that $\overline{s}_k=r_1$ and $\overline{s}_i = r_{k-i+1} - r_{k-i}$ if $i>1$.  This is the formula (\ref{E:duality 2}).
Finally, observe that the $j^{\rm th}$ group of overcrossing strands in  $\overline{\L}$, where $j=1,2,\dots,k$ intersects precisely $s_k+\cdots+s_{k-(j-1)}$ overcrossing strands of $\L$.  This gives the formula 
(\ref{E:duality 1}).  
\done

\begin{example}\label{Ex:dual vectors}\rm\quad We give some examples of dual Lorenz vectors.
\begin{enumerate}
\item $k=1$: \ $\langle r^s \rangle$ and $\langle s^r \rangle$ are dual vectors, so T$(r,s)=$T$(s,r)$.   
\item $k=2$: \ $\langle r_1^{s_1},\ r_2^{s_2}\rangle$  \ and \ $\langle s_2^{r_2-r_1},\ (s_1+s_2)^{r_1}\rangle$ \ are dual,\\ 
\hspace*{1.2cm} so T$((r_1,s_1),(r_2,s_2))\ = \ {\rm T}((s_2,r_2-r_1),(s_1+s_2,r_1))$.  
\item The example in Figure \ref{F:lorenz braid}: $\langle 2^4, 3^2, 6, 8^2 \rangle$ and $\langle 2^2, 3^3, 5, 9^2 \rangle$ are dual.
\item $\langle 2^2,\ 3^3,\ 4^2,\ 7,\ 9,\ 13^2\rangle $ \ and \ $\langle 2^4,\ 3^2,\ 4^3,\ 6,\ 9,\ 11^2\rangle$ are dual.
\end{enumerate}
\end{example}

\begin{remark}\rm
By defintion, Lorenz braids are positive.  However, T--links arise
naturally as a subset of generalized twisted torus links, which need
not be positive.  These are defined as in (\ref{E:T--braid}), except
we now allow $s_i\in\mathbb{Z}$; if $s_i<0$, the braid generators in
that syllable are negative.

Many of our results for T--links were obtained using the duality of
the Lorenz template.  Without positivity, there is no obvious duality,
but some of our results for T--links may still hold for generalized
twisted torus links.  

General twisted torus links are given by T$((r, s),(p, q))$ with
$p>r>0$.  If our duality holds, then $T((r, s),(p, q)) = T((q,
p-r),(q+s, r))$, which implies that $q+s>q>0$, hence $q,s>0$.
Therefore, duality as in Corollary \ref{C:duality} applies only to
positive twisted torus links.  

\begin{question}
Does another kind of duality apply to non-positive twisted torus links?
\end{question}
\end{remark}

\subsection{Upper bound for hyperbolic volume}
\label{S:hyperbolic}

Having the duality formulas of Corollary~\ref{C:duality} on hand, we
are ready to establish that the volume of hyperbolic Lorenz knots is
bounded by a constant that depends only on the size of the Lorenz
vector.  If $\mathcal{L}$ is obtained by Dehn surgery on a link $Y$,
then by Thurston's Dehn surgery theorem \c{Thurston}, the hyperbolic
volume of $\mathcal{L}$ is less than the hyperbolic volume of $Y$.
This theorem has many other implications that are easier to explore
using T--links.  For example, it follows that for any $1\leq i\leq k$,
there is a link $Y$ with an unknotted component whose volume is given
by
$$ \lim_{n\to\infty} {\rm Vol}({\rm T}((r_1,s_1),\ldots,(r_{i-1},s_{i-1}),(r_i,n\cdot r_i),(r_{i+1},s_{i+1}),\ldots,(r_k,s_k))).$$ 
Thurston's Dehn surgery theorem together with our results
shows that the volume is bounded for many infinite collections of
Lorenz links:
\begin{corollary}
\label{C:volume bound}
Let $N>0$.  Let $\mathcal{L}$ be a hyperbolic Lorenz link such
that its Lorenz vector has either $p\leq N$ or $d_p\leq N$; equivalently, its
T--braid has either $r_k\leq N$ or $(s_1+\ldots+s_k)\leq N$.  Then the
hyperbolic volume of $\mathcal{L}$ is bounded by a constant that
depends only on $N$.
\end{corollary}
\pf 
By Theorem \ref{T:Lorenz-braids and T-braids}, we can establish the
claim for T--links for which $r_k\leq N$ or $(s_1+\ldots+s_k)\leq N$.  Because
of the special form for T--braids in (\ref{E:T--braid}), we can
express the twists of $\mathcal{L}$ as Dehn surgeries on a nested
sequence of unknots, $\{(U_i,n_i)\}$, as in the proof of Corollary
\ref{C:applications2}.  Namely, for each $1\leq i\leq k$, we can
find some integers $n_i\geq 0,\ 0<a_i\leq r_i$ such that $s_i=n_i\cdot
r_i+a_i$.  Then, for all $1\leq i\leq k$, we perform a $1/n_i$--Dehn
surgery on $U_i$, augmented to
T$((r_1,s_1),\ldots,(r_{i-1},s_{i-1}),(r_i,a_i),(r_{i+1},s_{i+1}),\ldots,(r_k,s_k))$.

Therefore, if $r_k\leq N$, $\mathcal{L}$ is obtained by some Dehn
surgeries on a fixed finite collection of links.  For $\mathcal{L}$ such
that $s_1+\ldots+s_k\leq N$, by Corollary \ref{C:duality}, we consider the
dual T--link, with $\overline{r}_k\leq N$.  So every $\mathcal{L}$ is
obtained by Dehn surgeries on a fixed finite collection of links, which are
given by closed T--braids augmented with unknots.  The claim now
follows by Thurston's Dehn surgery theorem.  \done

\subsection{Polynomial invariants of Lorenz links}   
\label{S:Jones polynomial}

Polynomial invariants for certain infinite families of T--links are
known.  As another application of Theorem \ref{T:Lorenz-braids and
  T-braids}, we obtain the first such invariants for infinite families
of Lorenz links.
For the Jones polynomial, twisting formulas were given in Theorem 3.1
of \c{CK}.  Thus, the Jones polynomial of an infinite family of links
can be obtained from that of any one sufficiently twisted base case.

The Jones polynomials of Lorenz links are highly atypical.  The
polynomials are often sparse, nonzero coefficients are very small, and
the $L^1$-norm of coefficients is several orders of magnitude less
than for typical links with the same crossing number.  Mahler measure
is a natural measure on the space of polynomials for which these kinds
of polynomials are simplest.  Accordingly, the Mahler measure of Jones
polynomials of Lorenz links is unusually small, even when their span,
which is a lower bound for crossing number, is large.  

In \c{CK}, it was shown that the Mahler measure of the Jones
polynomial converges under twisting for any link:
Let $M(V_{\mathcal{L}})$ denote the Mahler measure of the Jones polynomial $V_{\mathcal L}$ of $\mathcal{L}$. 
For any $1\leq i\leq k$, there is a 2-variable polynomial $P$ such that
$$ \lim_{n\to\infty} M(V_{{\rm T}((r_1,s_1),\ldots,(r_{i-1},s_{i-1}),(r_i,n\cdot r_i),(r_{i+1},s_{i+1}),\ldots,(r_k,s_k))})  = M(P) $$

Thus, the atypical Jones polynomials of Lorenz links may be better
understood from the point of view of T--links.  For example, the
following result is similar to Corollary \ref{C:volume bound}:

\begin{corollary}
\label{C:Jones bound}
Let $N>0$.  Let $\mathcal{L}$ be a Lorenz link such that its
Lorenz vector has either $p\leq N$ or $d_p\leq N$; equivalently, its T--braid
has either $r_k\leq N$ or $s_1+\ldots+s_k\leq N$.  Then the Mahler measure of
the Jones polynomial of $\mathcal{L}$ is bounded by a constant that
depends only on $N$.
\end{corollary}
\pf The proof follows that of Corollary \ref{C:volume bound}, except
that at the end, the Dehn surgery theorem is replaced by the
$L^2$--bound for Mahler measure, as we explain below.

Let $L={\rm T}((r_1,a_1),\ldots,(r_k,a_k))$, as in the proof of
Corollary \ref{C:volume bound}.  For ${\mathbf n}=(n_1,\ldots,n_k)$,
construct $L_{\mathbf n}$ by $1/n_i$--Dehn surgeries on $U_i$, for $1\leq
i\leq k$.  By the proof of Corollary 2.3 of \c{CK}, there is a $(k+1)$--variable polynomial 
$P(t,x_1,\ldots,x_k)$ that depends only on $L$, such that $M(V_{L_\mathbf{n}}(t)) = M\left(P(t,t^{n_1},\ldots,t^{n_k})\right)$.
\footnote{
If we add $n$ full twists on $r$ strands of $L$, then the Kauffman
bracket polynomial $\kb{L_n}=A^{nr(r-1)}\,P(t,t^n)$, so $M(V_{L_n})=M(\kb{L_n})=M(P(t,t^n))$.
This is iterated for each twist site.}

If $||P||$ denotes the $L^2$--norm of coefficients of $P$, then $M(P)\leq ||P||$.  Therefore,
$$ M(V_{L_\mathbf{n}}(t)) = M\left(P(t,t^{n_1},\ldots,t^{n_k})\right) \leq ||P(t,t^{n_1},\ldots,t^{n_k})|| \leq ||P(t,x_1,\ldots,x_k)|| $$
So if $r_k\leq N$, $M(V_{\mathcal{L}})$ is bounded by $\max_j(||P_j||)$, for a fixed finite collection of polynomials $P_j$.
\done

Pierre Dehornoy has found many examples of distinct Lorenz knots with
the same Jones polynomial, with some pairs that have the same
Alexander polynomial as well \c{PDemail}.  For example,
$$ \langle 4, 4, 5, 7, 7, 7, 7, 7 \rangle \quad \text{and} \quad \langle 2, 3, 4, 5, 5, 6, 6, 6, 6, 6 \rangle $$
have a common Jones polynomial but different hyperbolic volume.  The
first knot above is also the knot $K7_{75}$ in the census of simplest
hyperbolic knots (see Section \ref{S:simplest hyperbolic}).  The Jones
polynomials of these and most other geometrically simple knots were
computed in \c{CKP}.

No general formula is known for Jones polynomials of Lorenz links,
even though calculations suggest that their Jones polynomials are very
special. We now give a statement that is true for all links that are
closed positive braids, and so in particular for all Lorenz links. Our
focus has been on the Jones polynomial, 
but it seems appropriate to also mention related results for the
Homflypt and Alexander polynomials P$_{\mathcal L}$ and $\Delta_{\mathcal
  L}$:
\begin{corollary}
  \label{C:Stoimenow} {\rm \cite{crom, Fi, Ka, St}}. Let $\mathcal L$
  be a link which is represented as a closed positive $n$--braid $L$.
  Let $c,\,\mu,\,u,\,g,$ be the number of crossings of $L$, the
  number of components, the unknotting number and the genus of $\mathcal L$.
  Then the following hold:
\begin{eqnarray*}
 2\min\deg({\rm V}_{\mathcal L}) & = & 2g+\mu-1 = 2u-\mu+1 = c-n+1\\
  & = & \max\deg(\Delta_{\mathcal L}) = \max\deg_z({\rm P}_{\mathcal L}) = \min\deg_v({\rm P}_{\mathcal L}) 
   \end{eqnarray*}
\end{corollary}
 
The Jones polynomials of twisted torus links T$((r_1,s_1),(r_2,s_2))$
are prime candidates for experiments because they are determined by
four integer parameters, i.e. $r_1,s_1,r_2,s_2$.  If we peek
ahead to Corollary~\ref{C:braid index of T--links} we will see that
the minimum braid index is a known function of these parameters.
Moreover, we know that any invariant, including
the Jones polynomial, must satisfy the duality of the defining
parameters.

\begin{question}
What is the Jones polynomial of T$((r_1,s_1),(r_2,s_2))$?
\end{question}

We turn briefly to the Alexander polynomial. By Theorem \ref{T:Lorenz-braids and T-braids}, we can find the
Alexander polynomial for a non-trivial infinite family of Lorenz links that are not torus links. We use the
fact that Morton \c{Morton}
computed the Alexander polynomial of T$((2,2m),(p,q))$,
$$ \Delta_{\rm T} = \frac{1-t}{(t^p-1)(t^q-1)}\left(1-(1-t)(1+t^2+\ldots+t^{2m-2})(t^a+t^b)-t^{pq+2m}\right) $$
Here, $a=pv$ and $b=(p-u)q$ where $0<u<p,\ 0<v<q$ and $uq\equiv -1\mod p$, $pv\equiv 1\mod q$.
By Theorem~\ref{T:Lorenz-braids and T-braids} this is the Alexander polynomial of the Lorenz link with defining vector $\langle 2^{2m}, p^q \rangle$. 

\subsection{Braid index formulae}
\label{S:Braid index}
By \c{FW}, the braid index $t$ of a Lorenz link is easily computed, one example
at a time, from the definition of the trip number $t$ that we gave in
(\ref{E:trip number}), but it is unclear how $t$ is related to the
parameters $\{(r_i,s_i), \ i=1,\dots,k\}$.  Our next application gives
a formula for the braid index which depends in a simple way on the
defining parameters. 

\begin{corollary}
\label{C:braid index of T--links}
Let $\mathcal{L}$ be the T--link T$((r_1,s_1),\ldots,(r_k,s_k))$.
Let $r_0= \overline{r}_0= 0$, so we can define using (\ref{E:duality 1}),
$$ i_0 = \min\{\ i\ |\ r_i \geq \overline{r}_{k-i}\} \quad {\rm and }\quad 
   j_0 = \min\{\ j\ |\ \overline{r}_j\geq r_{k-j}\} $$
Then the braid index of $\mathcal{L}$ is $ t(\mathcal{L}) = \min(r_{i_0},\, \overline{r}_{j_0})$.  

If $k=2$,  the braid index of any positive twisted torus link T$((r_1,s_1),(r_2,s_2))$ is given by:
$$ t = \begin{cases}
\min(s_2,\ r_2)     & {\rm if}\ r_1\leq s_2  \\
\min(s_1+s_2,\ r_1) & {\rm if}\ r_1\geq s_2
\end{cases}
$$

If $k=1$, i.e. torus links, our formula reduces to the well known fact that the braid index of T$(r,s)$ is $\min(r,s)$.
\end{corollary}
\pf
As a Lorenz link, $\mathcal{L}$ is defined by $\vec{d}_L = \langle r_1^{s_1},\ldots, r_k^{s_k}\rangle$.
Below, we use the notation in (\ref{E:1}) and (\ref{E:2}) with $r_{\mu_i} = d_i$,
so that the following are equivalent:
\begin{eqnarray*}
i + d_i & \geq & p \\
s_1 + \ldots + s_{\mu_i} + r_{\mu_i} & \geq & s_1 + \ldots + s_k \\
r_{\mu_i} & \geq & s_{\mu_i + 1} + \ldots + s_k \\
r_{\mu_i} & \geq & \overline{r}_{k-\mu_i}
\end{eqnarray*}
Therefore, $\displaystyle i_0 = \min\{\ \mu_i\ |\ i + d_i \geq p \}$ and 
$\displaystyle j_0 = \min\{\ \mu_j\ |\ j + \overline{d}_j \geq \overline{p}=d_p \} $.

Since displacements correspond to intersecting strands, the $i$-th
overcrossing strand crosses $d_i$ undercrossing strands.  Thus, by
(\ref{E:trip number}), $t$ is the number of LR-strands, which equals
the number of RL-strands.  We now consider two cases.

{\bf Case 1.}\quad 
There exists $i_*$ such that $i_* + d_{i_*} = p$.

The left strand $\alpha$ starting at $i_*$ with endpoint $p$ is the
last LL-strand, so it does not intersect any RR-strands.  The
equality implies that all RL strands intersect $\alpha$, so $d_{i_*}=|$RL$|$.  
For all $i<i_*$, $i+d_i < p$, so $i_0 \geq \mu_{i_*}$.
For all $i\geq i_*$, $r_{\mu_i}\geq |$RL$|=r_{\mu_{i_*}}$, so $i_0 = \mu_{i_*}$. Therefore,
$$ t=|{\rm RL}|=d_{i_*}=r_{i_0} $$

{\bf Case 2.}\quad 
There does not exist $i_*$ such that $i_* + d_{i_*} = p$.

There exists a right strand $\gamma$ with endpoint $p$, which is the
first RL-strand.  Because its endpoint is $p$, $\gamma$ intersects
all LR-strands and no LL strands.  In the dual Lorenz link $\overline{\mathcal{L}}$,
if $\gamma$ starts at $j_*$ then $\overline{d}_{j_*}=|$LR$|$.
By duality, the endpoint of $\gamma$ is $\overline{p}+1$.
If another strand $\gamma'$ is parallel to $\gamma$ with endpoint $\overline{p}$ then 
both strands are in the same packet, so $t=\overline{d}_{j_*}=\overline{r}_{j_0}$ by Case 1 applied to $\overline{\mathcal{L}}$.
Otherwise, for all $j<j_*$, $\overline{d}_j< \overline{d}_{j_*}$ so $j + \overline{d}_j< \overline{p}$, hence $j_0\geq \mu_{j_*}$. 
For all $j\geq j_*$, $\overline{r_j}\geq |LR|=\overline{r}_{\mu_{j_*}}$, so $j_0= \mu_{j_*}$.
Therefore, $$ t=|{\rm LR}|=\overline{d}_{j_*}=\overline{r}_{j_0} $$ 

In both cases, $r_{i_0},\overline{r}_{j_0} \geq t$, so  $t=\min(r_{i_0}, \overline{r}_{j_0})$.

Now let's specialize to the case $k=2$. 
If $r_1 \geq s_2 = \overline{r}_1$ then either case below occurs:
\begin{enumerate}
\item[$i$.] $r_1 \geq \overline{r}_2 > \overline{r}_1 \quad \Rightarrow\quad t = \overline{r}_2 =  s_1+s_2$
\item[$ii$.] $\overline{r}_2 \geq r_1 \geq \overline{r}_1 \quad \Rightarrow \quad t = r_1 $
\end{enumerate}
If $\overline{r}_1= s_2  \geq r_1$  then either case below occurs:
\begin{enumerate}
\item[$iii$.] $r_2 \geq \overline{r}_1 \geq r_1  \quad \Rightarrow \quad t = \overline{r}_1 =  s_2$
\item[$iv$.] $\overline{r}_1 \geq r_2 > r_1 \quad \Rightarrow \quad t = r_2 $
\end{enumerate}

When $k=1$, the Lorenz braid defined by $d_L = \langle r^s \rangle$
represents the torus link T$(r,s)$, which is the closure of the
$r$-braid $(\sigma_1\dots\sigma_{r-1})^s$.  The dual Lorenz braid
$\langle s^r \rangle$ represents the same torus link T$(s,r)$, which
is the closure of the $s$-braid $(\sigma_1\dots\sigma_{s-1})^r$.  As
is well known, the braid index of a torus link is $\min(r,s)$, which
agrees with Corollary~\ref{C:braid index of T--links}.
\done

\section{Minimal braid index representations} 
\label{S:Minimal braid index representations}

We have proved that there are different closed braid representations
of a Lorenz link $\mathcal{L}$, with braid indices: $p+d_p,\ d_p,\ p$,
and $t$.  The representation of braid index $p+d_p$ is the Lorenz
braid defined by our vector $d_L = \langle r_1^{s_1},\dots, r_k^{s_k}
\rangle$.  The representation of braid index $d_p=r_k$ was given in
Theorem~\ref{T:Lorenz-braids and T-braids}, and its dual $p$--braid in
Corollary~\ref{C:duality}.  In this section, we use
Lemma~\ref{L:Markov sequence} and some of the things we have learned
along the way, to establish another correspondence, 
this time between Lorenz braid representations of $\mathcal{L}$ (hence
also T--braid representations) and special $t$--braid representations,
where $t$ is the minimal braid index of $\mathcal{L}$.

Let $\mathbb{L}$ be a Lorenz braid defined by $\vec{d}_L=\langle
d_1,d_2,\dots,d_p \rangle$, representing the Lorenz link
$\mathcal{L}$.  As discussed earlier, the strands in $\mathbb{L}$ divide into
strands of type LL, LR,RL and RR, where strand $j$ has type LL if and
only if $1\leq j\leq p-t$.  By duality, strand $j$ has type RR if and
only if strand $\overline{j}$ has type LL with respect to
$\overline{\mathbb{L}}$; i.e., $1\leq \overline{j}\leq\overline{p}-t
=d_p-t$.  We define
\begin{eqnarray}
\label{E:ni}
n_i &=& \# \{ {\rm strand}\ j \in {\rm LL \quad such \ that} \quad d_j = i+1\}, \\
\label{E:mi}
m_i &=& \# \{ {\rm strand}\ j \in {\rm RR \quad  such \ that} \quad  \overline{d}_{\overline{j}} = i+1\}.
\end{eqnarray}

Let $\vec{n} = (n_1, \dots,n_{t-1})$, $\vec{m} = (m_1,\dots,m_{t-1})$,
where each $n_i,\; m_j \geq 0$.  The conditions in (\ref{E:2}) are
automatically satisfied for any $\vec{n},\ \vec{m}$ with non-zero
entries.  In Example \ref{Ex:favorite}, for which $t=3$, we get
$n_1=4,\ n_2=2,\ m_1=2,\ m_2=3$, which is immediate from
Figure~\ref{F:lorenz braid}.

The triple $(t,\vec{n}, \vec{m})$, $t\geq 2$, defines the following
$t$--braid representation of $\mathcal{L}$, where $t$
is the braid index of $\mathcal{L}$:
\begin{eqnarray}
\label{E;t-braid formula}
\nonumber
 \mathbb{M} & = &  (\s 1\cdots \s {t-1})^t \ \  \prod_{i=1}^{t-1}\, (\s 1 \cdots \s i)^{n_i} \  
\prod_{i=t-1}^{1}\,  (\s {t-1} \cdots \s i)^{m_{t-i}} \\
\nonumber
 & = & [1,t]^t \ \ \prod_{i=1}^{t-1}\,  [1,i+1]^{n_i} \ \prod_{i=t-1}^1\,  [t,i]^{m_{t-i}} \\
& = &  [1,t]^t \ \ \prod_{i=1}^{t-1}\,  [1,i+1]^{n_i} \ \ \prod_{i=1}^{t-1}\  [t, t-i]^{m_i}
\end{eqnarray}
This was proved in Proposition 5.6 of \cite{BW}, with a small but very
confusing typo\footnote{In Proposition 5.6 of \cite{BW}, the
  superscript in the product on the right should have been $m_{t-i}$,
  not $m_i$.}
corrected.  The proof of Proposition 5.6 of \cite{BW} is correct, but the formula there is not.
The following theorem is a strengthening of Proposition 5.6 of
\cite{BW}, and is comparable to Theorem~\ref{T:Lorenz-braids and T-braids}: It sets up a 
correspondence between Lorenz braid representations of
$\mathcal{L}$ and special $t$--braid representations.

\begin{theorem}\label{T:Lorenz-braids and t-braids}
  There is a one-to-one correspondence between Lorenz braids
  $\mathbb{L}\in B_{p+d_p}$, with defining vector $\vec{d}_L$ as in
  (\ref{E:1}), and triples $(t,\vec{n},\vec{m})$, which determine a
  $t$-braid $\mathbb{M}$.
\end{theorem}

{\bf Caution:}  Distinct Lorenz braid representations of $\mathcal L$, and their
  corresponding distinct triples $(t,\vec{n},\vec{m})$, may determine
  the same $t$-braid representation $\mathbb M$ of $\mathcal L$.  For
  example, when $t=2$ the only possibility is $\vec{n} = (n_1), \ \ \vec{m} = (m_1)$,
  where $n_1$ and $m_1$ are both positive, so that $\mathbb M = \sigma_1^{2 + n_1 + m_1}$.  
  Other partitions of  $n_1+m_1$ will give other triples $(t,\vec{n},\vec{m})$ but the same
  2-braid.

\pf The reader is referred to \cite{BW} for the proof that $\mathbb L$
determines $\mathbb M$. We will prove the converse.

We first prove that $(t,\vec{n},\vec{m})$ determines the subvector
$\vec{d}_{LL} \subset \vec{d}_L$ consisting of all $d_i$ such that
$i\leq |$LL$|$.  In the proof of Lemma~\ref{L:Markov sequence}, we
uncoiled the LL--braid to construct the equivalent $t$--strand braid
$\mathbb{X}$, given in (\ref{E:expression for X}).  Namely, we traded
each braid strand in LL, together with its associated loop around the
axis, for an arc corresponding to one of the sequences $(\s 1 \s 2
\cdots \s {d_j-1})$ in the braid word $\mathbb{X}$.  The definition of
$n_i$ implies that
$$ \mathbb{X} = \prod_{j=1}^{p-t} (\s 1 \s 2 \cdots \s {d_j-1}) = \prod_{i=1}^{t-1} (\s 1 \s 2 \cdots \s i)^{n_i} $$
This is a subword of $\mathbb{M}$.  Going the other way, each subword
$(\s 1 \s 2 \cdots \s i)^{n_i} \in\mathbb{M}$ must have come from a
group of $n_i$ parallel strands in LL.  Since the LL braid is made up
entirely from groups of parallel strands, it follows that $|$LL$| = n_1 +
\cdots + n_{t-1}$.

Using the now-familar trick of passing from $\mathbb{L}$ to
$\overline{\mathbb{L}}$, it follows that $|$RR$| = m_1+\cdots +
m_{t-1}$.  Note also that, since $p = |$LL$| + t$ and $d_p = |$RR$| +
t$, it follows that the braid index of $\mathbb{L}$ is $$p + d_p = 2t
+ \sum_{i=1}^{t-1}(n_i + m_i). $$

Therefore, $(t,\vec{n},\vec{m})$ determines (i) the braid index of
$\mathbb{L}$, (ii) the number $|$LL$|$ of strands in the LL braid and
(iii) the subvector $\vec{d}_{LL} \subset \vec{d}_L$.  By duality,
$(t,\vec{m},\vec{n})$ then also determines $\vec{d}_{RR}$.  Next,
notice that the only strands of $\mathbb{L}$ which have endpoints in R
are type RR and LR, and from this it follows that all endpoint
positions in R which are not occupied by strands of type RR must be
occupied by the strands of type LR.  Moreover, the endpoints of the LR
strands are completely determined because there are no crossings
between pairs of strands of type LR.  Since we already know the vector
$\vec{d}_{LL}$, it follows that the vector $\vec{d}_L$ is completely
determined.  Likewise, $\vec{d}_R$ is determined, hence $\mathbb{L}$
is completely determined by $(t,\vec{n},\vec{m})$.
\done

\begin{remark} \rm 
  Using Theorem \ref{T:Lorenz-braids and t-braids}, we get a second
  proof of Corollary \ref{C:finite1}.  By
  Corollary~\ref{C:applications1}~$(ii)$, for fixed braid index, the
  letter length of any braid representation is a topological invariant
  of $\mathcal{L}$.  Let $t$ be the trip number of $\mathcal{L}$.
  Since only finitely many positive words have given letter length,
  there are finitely many $t$--braid representations of $\mathcal{L}$
  of the form (\ref{E;t-braid formula}).  By
  Theorem~\ref{T:Lorenz-braids and t-braids}, $\mathcal{L}$ has
  finitely many Lorenz braid representations of the form (\ref{E:1});
  i.e., up to trivial stabilizations.
\end{remark}

\begin{remark}\label{R:duality and t-braids} \rm
  Corollary~\ref{C:duality} results in a duality for $t$--braids,
  given by conjugation by the half-twist $\Delta$, which sends every
  $\sigma_i$ to $\sigma_{t-i}$.  For every $t$--braid as in
  (\ref{E;t-braid formula}), we get another braid in the same
  conjugacy class and which has the special form given in
  (\ref{E;t-braid formula}).  To see this, note that conjugation by
  $\Delta$ sends
$$\mathbb{M} = [1,t]^t \prod_{i=1}^{t-1} [1,i+1]^{n_i}\prod_{i=1}^{t-1} [t,t-i]^{m_{i}}$$ to 
$$\Delta\mathbb{M}\Delta^{-1} = [t,1]^t \prod_{i=1}^{t-1} [t,t-i]^{n_{i}}\prod_{i=1}^{t-1} [1,i+1]^{m_i} \approx [1,t]^t \prod_{i=1}^{t-1}[1,i+1]^{m_i} \prod_{i=1}^{t-1} [t,t-i]^{n_{i}}$$ 
where $\approx$ means after cyclic permutation.  We
use the fact that $\Delta^2=[1,t]^t = [t,1]^t$ is in the center of $B_t$.
\end{remark}

Our experimental data suggests that this is a general phenomenon:
\begin{conjecture} \label{C:unique conjugacy classes} If 
  a Lorenz link $\mathcal{L}$ has representations $\mathbb{M}_1, \
  \mathbb{M}_2\in B_t,$ where $t$ is the trip number of $\mathcal{L}$, then
  $\mathbb{M}_1,\mathbb{M}_2$ are in the same conjugacy class in $B_t$.
\end{conjecture} 

With regard to this conjecture, Corollary~\ref{C:applications2}
provides many examples of interesting conjugacy between the $t$--braid
representations of $\mathbb{L}_1$ and $\mathbb{L}_2$.  In general,
links that are closed positive braids need not have unique conjugacy
classes of minimum braid representations, but the known examples that
might contradict Conjecture~\ref{C:unique conjugacy classes} cannot be
Lorenz links. For example, composite links have minimum closed braid
representations that admit exchange moves, leading to infinitely many
conjugacy classes of minimum braid index representations, but Lorenz
links are prime \cite{W}.  Also, links that are closed 3-braids and
admit positive flypes have non-unique conjugacy classes, but the
Lorenz links of trip number 3 have been studied \cite{Bedient}, and
they do not include any closed positive 3-braids that admit positive
flypes.

There are very few families of links for which we know, precisely,
minimum braid index representatives, the most obvious being the unknot
itself.  In \cite{BW} Lorenz braids whose closures define the unknot
were delineated precisely.  The question of which Lorenz braids
determine torus links is more complicated, but is a natural next step.
A pair of positive integers $p,q$ suffice to determine the type of any
torus link, but looking at the class of all Lorenz links, it is
difficult to determine which ones are torus links.  With the help of
Theorem \ref{T:Lorenz-braids and t-braids}, we are able to make a
contribution to that problem:

\begin{corollary}
\label{C:algorithm}
Let $\mathcal L$ be a Lorenz link with trip number $t$.  Let $\mathbb
M$ be a $t$-braid representative of $\mathcal L$, as given in
(\ref{E;t-braid formula}).  Then there is an algorithm of complexity
$O(|\mathbb M|^2 t^3\log t)$ that determines whether or not $\mathcal
L$ is a torus link.
\end{corollary}

\pf 
By \cite{FW}, $t$ is the braid index of
$\mathcal L$. By a theorem of Schubert \cite{Sc}, we know that if
$\mathcal L$ is a torus link, then it has a minimum braid index
representative in $B_t$ of the form $[1,t]^q$ for some $q\geq t$.
Also, by a different result in \cite{Sc}, any closed $t$-braid that
represents $\mathcal L$ must be conjugate to $[1,t]^q$.  Our first
question is: if $\mathcal L$ is a torus link, what is the integer $q$?
Since $\mathbb M$ and $[1,t]^q$ are both positive braids, they must
have the same letter length.  From this it follows that $\mathcal L$
cannot be a torus link unless $|\mathbb M| = (t-1)(q)$ . Therefore
$q = |\mathbb M|/(t-1)$.

We now give an algorithm to decide whether the $t$-braids $\mathbb
L_t$ and $[1,t]^q$ are conjugate in $B_t$.  Changing our viewpoint, we
now regard the braid group $B_t$ as the mapping class group of the
unit disc $D^2$ minus $t$ points, where admissible maps fix
$\partial D^2$.  See \cite{BB}, for example, for a proof that
this mapping class group is isomorphic to Artin's braid group $B_t$.
Let $\delta$ be the $t$-braid $[1,t]$.  If the points which are
deleted from the unit disc $D^2$ are arranged symmetrically
around the circle of radius 1/2, then $\delta$ may be seen as a
rotation of angle $2\pi/t$ about the origin, with the boundary of
$D^2$ held fixed.  Such a braid has the Thurston-Nielson type
of a {\it periodic} braid of period $t$.  By the results in \cite{GM}
we know that periodic braids have unique roots. Therefore it suffices
to prove that $(\mathbb M)^t$ is conjugate to $\delta^{qt}$.
Observe that $\delta^t$ generates the center of $B_t$, so
that $\delta^{tq}$ is in the center, and from this it follows that it
suffices to prove that $\mathbb M^t$ and $\delta^{tq}$ represent the
same element of $B_t$.  This trick reduces the conjugacy problem to
the word problem.

There is a solution to the word problem in $B_t$ that was discovered
simultaneously by El-Rifai--Morton and by Thurston which has the
property: if $X$ is an element of $B_t$ which has letter length $|X|$
then its {\it left-greedy normal form} can be computed in time
$O(|X|^2 t \log t)$.  In our case the word length of $\mathbb M^t$
is $(t)(|\mathbb M|)$, therefore the problem can be solved in time
$O(|\mathbb M|^2 t^3 \log t)$, as claimed.  
\done

\begin{remark}\rm
\label{R:satellite}
El-Rifai \c{ER} classified all ways in which a Lorenz knot can be
presented as a satellite of a Lorenz knot.  He showed that only
parallel cables with possible twists can occur.  These results
generalize Theorems 6.2 and 6.5 of \c{BW}.   
\begin{question}
\label{Q:satellite1}
  Is there an efficient algorithm, along the lines of
  Corollary~\ref{C:algorithm}, to recognize
  when a Lorenz knot is a satellite of a Lorenz knot?
\end{question}

In relation to the above, a very interesting open problem was posed in
\cite{ER}: 
\begin{question}
\label{Q:satellite2}
  Can a Lorenz knot be a satellite of a non-Lorenz knot?
\end{question}
Noting the method of proof in \cite{W} that Lorenz knots are
prime, one suspects that the fact that every Lorenz knot embeds on the
Lorenz template implies that the answer is `no'.  This is an
interesting question because one would like very much to know how to
separate the hyperbolic Lorenz knots and links from those which are
not hyperbolic.

In this regard, we note that the Lorenz braids that determine the
unknot were completely characterized in \cite{BW}.  It seems to be
much more difficult to decide: 
\begin{question}
  Which Lorenz braids close to torus links?
\end{question}
We have partial results on this problem, but have not
found a satisfactory general answer.
\end{remark}

\section{Lorenz data for the simplest hyperbolic knots}
\label{S:simplest hyperbolic}

In the table below, we list 107 simplest hyperbolic knots (see \c{CDW,
  CKP}) that are Lorenz, and five that are possibly Lorenz; the rest are
not Lorenz.  The symbol $\mathbf{k}n_{m}$ means the $m$th knot in the
census of hyperbolic knots whose complement can be constructed from
no less than $n$ ideal tetrahedra.

The 107 identified Lorenz braids in the table were proved to be
isometric to the corresponding census knots using SnapPea to verify
the isometry.  Many had already been identified as positive twisted
torus knots in \c{CDW, CKP}.

The 89 simplest hyperbolic knots that are not listed in the table are
not Lorenz.  For many, their Jones polynomials from \c{CKP} failed to
satisfy Corollary~\ref{C:Stoimenow}.  For the others, we used the
following method:  

Pierre Dehornoy computed all Lorenz braids up to 49 crossings that
close to a knot, and Slavik Jablan eliminated duplications from this
list, which contains 14,312 distinct non-alternating Lorenz knots up
to 49 crossings.  If $c$ is the crossing number of the Lorenz braid
and $g$ is the genus of the knot, then by (\ref{E:c4g}), we know that
$c\leq 4g$.  So any Lorenz knot with $g\leq 12$ has a Lorenz braid
representation with $c\leq 48$.  Therefore, any knot with $g\leq 12$
that is missing from the Dehornoy--Jablan list cannot be Lorenz.

Knots with 16 or fewer crossings are classified, and their invariants
are accessible using Knotscape.  For these knots, if the minimal and
maximal degrees of their Jones polynomials have the same sign, we
verified that the smaller absolute value of the two is less than 12.
It follows that $g\leq 12$ for any of these knots that satisfy
Corollary~\ref{C:Stoimenow}.  Jablan provided us with the
following Knotscape knots in the Dehornoy--Jablan list, which is,
therefore, the complete classification of Lorenz knots up to 16 crossings:
\begin{eqnarray*}
 3_{1},\, 5_{1},\, 7_{1},\, 9_{1},\, 8_{19},\, 10_{124},\, 11a_{367},\, 12n_{242},\, 12n_{725},\, 13a_{4878},\, 14n_{6022},\, 14n_{21324} \\
14n_{21881},\, 15n_{41185},\, 15a_{85263},\, 16n_{184868},\, 16n_{771177},\, 16n_{783154},\, 16n_{996934} 
\end{eqnarray*}

In addition, Jablan verified for us that $\mathbf{k}6_{31}$,
$\mathbf{k}7_{81}$, $\mathbf{k}7_{83}$, $\mathbf{k}7_{106}$,
$\mathbf{k}7_{113}$, $\mathbf{k}7_{118}$, $\mathbf{k}7_{124}$, and
$\mathbf{k}7_{119}$ are not on the Dehornoy--Jablan list.  The
remaining four simplest hyperbolic knots, indicated by a ``?'' in the
table below, have diagrams with more than 49 crossings, which cannot
be handled by this computer program.  Except for $\mathbf{k}7_{119}$,
the knots listed have Jones polynomials (see \c{CKP}) that imply
$g\leq 12$ if they satisfy Corollary~\ref{C:Stoimenow}.  Although
$\mathbf{k}7_{119}$ has a diagram with 33 crossings, $g=15$, so we
cannot be certain that it does not have a Lorenz braid with $50\leq
c\leq 60$.

The following formulas, which follow from results earlier in this paper, provide additional
information that can be obtained using the Lorenz braids in our final table on the next page.
Let $\mathcal{L}$ be any Lorenz link given by $\vec{d}_L = \langle
d_1,\dots,d_p\rangle$, as in (\ref{E:1}).  Let $S=\sum_{i=1}^p d_i$,
and let $t$ be its trip number.  The crossing
numbers and braid indices of the Lorenz braid $\L$, the T--braid $\TT$,
the dual T--braid $\TT'$, and the minimal braid index $t$--braid ${\mathbb M}$ are as follows:

\ 

\renewcommand{\arraystretch}{1.15}
\centerline{
\begin{tabular}{|c|c|c|c|c|}
\hline 
&$\mathbb L$ & $\mathbb T$ & $\mathbb T'$ & $\mathbb M$\\ 
\hline
crossing number & $S$     & $S-p$ & $S-d_p$ & $S+t-p-d_p$\\ 
\hline
braid index & $p+d_p$ & $d_p$ & $p$     & $ t $\\
\hline
\end{tabular}}

\ 

The braid crossing numbers of the braids $\mathbb M$ for some Lorenz
knots in the census turn out to be surprisingly high.   In fact, the crossing number of the minimum index braid in equation (\ref{E;t-braid formula}) is the minimal crossing number of the Lorenz link, by Proposition 7.4 of \cite{Murasugi}.

On the next page, we give the table of Lorenz knots that are in the census of hyperbolic knots whose complements 
can be constructed from seven or fewer ideal tetrahedra.

\newpage
\thispagestyle{empty}
\renewcommand{\arraystretch}{1.15}
\centerline{
\begin{tabular}{|c |c | c | c | c |c | c | c | }
\hline
Knot & Lorenz vector & \ \ & Knot & Lorenz vector & \ \  & Knot & Lorenz vector \\
\hline
\hline
$\mathbf{k}3_{1}$ & $\kb{2^{2},\, 3^{4}}$ &&$\mathbf{k}7_{3}$ & $\kb{2^2,\, 5^{16}}$ 
&&$\mathbf{k}7_{55}$ & $\kb{5^{2},\, 7^9}$   \\ \hline
$\mathbf{k}4_{3} $ & $\kb{2^{2},\, 3^{8}}$ &&$\mathbf{k}7_{4}$ & $\kb{3^3,\, 5^{17}}$
&&$\mathbf{k}7_{56}$ & ? \\\hline
$\mathbf{k}4_{4} $ & $\kb{2^{2},\, 4^{7}}$ &&$\mathbf{k}7_{5}$ & $\kb{2^4,\, 3^{16}}$
&&$\mathbf{k}7_{57}$ & $\kb{3^{6},\, 4^{13}}$ \\\hline
$\mathbf{k}5_{1} $ & $\kb{2^{2},\, 5^{6}}$ &&$\mathbf{k}7_{6}$ & $\kb{3^2,\, 7^{16}}$
&&$\mathbf{k}7_{58}$ & $\kb{4^{4},\, 7^{10}}$ \\ \hline
$\mathbf{k}5_{4} $ & $\kb{2^{2},\, 5^{8}}$ &&$\mathbf{k}7_{7}$ & $\kb{2^2,\, 7^{18}}$
&&$\mathbf{k}7_{59}$ & $\kb{5^{4},\, 9^{7}}$ \\\hline
$\mathbf{k}5_{5} $ & $\kb{2^{2},\, 3^{11}}$ &&$\mathbf{k}7_{8}$ & $\kb{2^2,\, 6^{7}}$
&&$\mathbf{k}7_{60}$ & $\kb{7^{7},\, 8^{5}}$  \\\hline
$\mathbf{k}5_{6} $ & $\kb{3^{3},\, 5^{6}}$  &&$\mathbf{k}7_{9}$ & $\kb{10^2,\, 11^{4}}$
&&$\mathbf{k}7_{61}$ & $\kb{7^{10},\, 3^{2}}$   \\ \hline
$\mathbf{k}5_{7} $ & $\kb{2^{2},\, 5^{7}}$  &&$\mathbf{k}7_{12}$ & $\kb{3^2,\, 8^{7}}$
&&$\mathbf{k}7_{62}$ & $\kb{8^{8},\, 10^{3}}$  \\ \hline
$\mathbf{k}5_{10} $ & $\kb{4^{2},\, 5^{4}}$ &&$\mathbf{k}7_{13}$ & $\kb{2^2,\, 7^{12}}$
&&$\mathbf{k}7_{63}$ & $\kb{4^{2},\, 6^{11}}$   \\ \hline
$\mathbf{k}5_{11} $ & $\kb{2^{6},\, 3^{4}}$ &&$\mathbf{k}7_{14}$ & $\kb{2^2,\, 8^{11}}$
&&$\mathbf{k}7_{64}$ & $\kb{4^{4},\, 5^{17}}$   \\ \hline
$\mathbf{k}5_{14} $ & $\kb{3^{2},\, 4^{7}}$ &&$\mathbf{k}7_{15}$ & $\kb{10^4,\, 11^{4}}$
&&$\mathbf{k}7_{66}$ & $\kb{5^{3},\, 6^{11}}$    \\ \hline
$\mathbf{k}5_{15} $ & $\kb{4^{4},\, 7^{3}}$ &&$\mathbf{k}7_{16}$ & $\kb{5^3,\, 8^{11}}$
&&$\mathbf{k}7_{67}$ & $\kb{8^{4},\, 9^{10}}$   \\ \hline                    
$\mathbf{k}5_{16} $ & $\kb{5^{5},\, 7^{3}}$ &&$\mathbf{k}7_{17}$ & $\kb{3^2,\, 8^{13}}$      
&&$\mathbf{k}7_{68}$ & $\kb{2^{6},\, 3^{10}}$    \\ \hline
$\mathbf{k}5_{17} $ & $\kb{4^{4},\, 5^{7}}$ &&$\mathbf{k}7_{20}$ & $\kb{5^2,\, 6^{5}}$     
&&$\mathbf{k}7_{69}$ & $\kb{7^{4},\, 8^{9}}$    \\ \hline
$\mathbf{k}5_{18} $ & $\kb{4^{8},\, 5^{3}}$ &&$\mathbf{k}7_{21}$ & $\kb{2^2,\, 7^{9}}$      
&&$\mathbf{k}7_{71}$ & $\kb{6^{6},\, 7^{10}}$   \\ \hline
$\mathbf{k}6_{3} $ & $\kb{2^{2},\, 5^{11}}$ &&$\mathbf{k}7_{22}$ & $\kb{3^3,\, 5^{16}}$     
&&$\mathbf{k}7_{73}$ & $\kb{8^{6},\, 9^{8}}$   \\ \hline
$\mathbf{k}6_{4} $ & $\kb{3^{3},\, 5^{12}}$ &&$\mathbf{k}7_{23}$ & $\kb{5^{15},\, 7^{2}}$      
&&$\mathbf{k}7_{75}$ & $\kb{3^{3},\, 5^{8}}$   \\ \hline
$\mathbf{k}6_{5} $ & $\kb{2^{2},\, 3^{14}}$ &&$\mathbf{k}7_{27}$ & $\kb{3^{2},\, 4^{15}}$
&&$\mathbf{k}7_{76}$ & $\kb{6^{12},\, 7^{4}}$   \\ \hline
$\mathbf{k}6_{6} $ & $\kb{3^{2},\, 7^{9}}$ &&$\mathbf{k}7_{28}$ & $\kb{4^{8},\, 5^{6}}$ 
&& $\mathbf{k}7_{78}$ &  $\kb{6^{12},\, 7^{5}}$    \\ \hline
$\mathbf{k}6_{7} $ & $\kb{2^{2},\, 7^{11}}$ &&$\mathbf{k}7_{29}$ & $\kb{7^{2},\, 9^{11}}$
&&$\mathbf{k}7_{79}$ & $\kb{3^{2},\, 9^{13}}$   \\ \hline
$\mathbf{k}6_{11} $ & $\kb{2^{2},\, 5^{9}}$  &&$\mathbf{k}7_{30}$ & $\kb{2^{2},\, 4^{19}}$
&&$\mathbf{k}7_{82}$ & $\kb{4^{4},\, 5^{13}}$    \\ \hline
$\mathbf{k}6_{12} $ & $\kb{3^{3},\, 5^{11}}$ &&$\mathbf{k}7_{31}$ & $\kb{5^{3},\, 9^{13}}$
&&$\mathbf{k}7_{87}$ & $\kb{5^{5},\, 7^{4}}$   \\ \hline
$\mathbf{k}6_{13} $ & $\kb{2^{2},\, 5^{12}}$ &&$\mathbf{k}7_{32}$ & $\kb{8^{4},\, 11^{5}}$
&&$\mathbf{k}7_{88}$ & $\kb{5^{15},\, 7^{3}}$   \\ \hline
$\mathbf{k}6_{14} $ & $\kb{3^{2},\, 4^{11}}$ &&$\mathbf{k}7_{33}$ & $\kb{6^{2},\, 7^{10}}$
&&$\mathbf{k}7_{90}$ & $\kb{2^{6},\, 3^{8}}$   \\ \hline
$\mathbf{k}6_{15} $ & $\kb{6^{6},\, 7^{4}}$ &&$\mathbf{k}7_{34}$ & $\kb{4^{2},\, 5^{9}}$ 
&&$\mathbf{k}7_{99}$ & $\kb{2^{2},\, 3^{2},\, 5^{2}}$   \\ \hline
$\mathbf{k}6_{16} $ & $\kb{7^{7},\, 8^{3}}$  &&$\mathbf{k}7_{35}$ & $\kb{7^{7},\, 9^{4}}$&&
$\mathbf{k}7_{101}$ & ?    \\ \hline
$\mathbf{k}6_{17} $ & $\kb{2^{2},\, 4^{15}}$ &&$\mathbf{k}7_{36}$ & $\kb{6^{6},\, 7^{9}}$
&&$\mathbf{k}7_{102}$ & $\kb{2^{4},\, 4^{5}}$   \\ \hline
$\mathbf{k}6_{18} $ & $\kb{6^{6},\, 7^{5}}$ &&$\mathbf{k}7_{37}$ & $\kb{2^{2},\, 5^{14}}$
&&$\mathbf{k}7_{109}$ &  ?  \\ \hline
$\mathbf{k}6_{19} $ & $\kb{2^{6},\, 3^{5}}$ &&$\mathbf{k}7_{38}$ & $\kb{2^{2},\, 8^{13}}$
&&$\mathbf{k}7_{110}$ & $\kb{2^2,\,3^5,\,5^4}$   \\ \hline
$\mathbf{k}6_{21} $ & $\kb{5^{4},\, 8^{7}}$ &&$\mathbf{k}7_{39}$ & $\kb{3^{3},\, 11^{2}}$
&& $\mathbf{k}7_{111}$ & $\kb{3^{3},\, 7^{10}}$  \\ \hline
$\mathbf{k}6_{25} $ & $\kb{4^{6},\, 5^{4}}$ && $\mathbf{k}7_{42}$ & $\kb{7^{14},\, 8^{3}}$
&& $\mathbf{k}7_{112}$ & $\kb{2^1,\,5^4,\,8^6}$   \\ \hline
$\mathbf{k}6_{27} $ & $\kb{3^{3},\, 4^{11}}$&&$\mathbf{k}7_{43}$ & $\kb{6^{2},\, 7^{16}}$
&& $\mathbf{k}7_{115}$ & $\kb{3^{5},\, 4^{7}}$  \\ \hline
$\mathbf{k}6_{29} $ & $\kb{3^{6},\, 7^{2}}$  &&$\mathbf{k}7_{47}$ & $\kb{4^{4},\, 6^{7}}$ 
&& $\mathbf{k}7_{119}$ &   ?   \\ \hline
$\mathbf{k}6_{30} $ & $\kb{4^{4},\, 5^{12}}$ &&$\mathbf{k}7_{48}$ & ?                                
&& $\mathbf{k}7_{122}$  & $\kb{4^2,\,5^3,\,7^3}$ \\ \hline
$\mathbf{k}6_{32} $ &  $\kb{4^{4},\, 5^{8}}$    &&$\mathbf{k}7_{50}$ & $\kb{6^{4},\, 7^{5}}$
&& $\mathbf{k}7_{123}$ & $\kb{6^{12},\, 8^{5}}$ \\ \hline
$\mathbf{k}6_{35} $ & $\kb{6^{6},\, 5^{8}}$   &&$\mathbf{k}7_{51}$ & $\kb{5^{2},\, 7^{8}}$
&& $\mathbf{k}7_{126}$ &  $\kb{4^{8},\, 8^{3}}$  \\ \hline
$\mathbf{k}6_{36} $ &  $\kb{5^{10},\, 7^{3}}$  &&$\mathbf{k}7_{52}$ & $\kb{5^{2},\, 7^{6}}$
&& &   \\  \hline
$\mathbf{k}6_{39} $ &$\kb{4^{4},\, 8^{3}}$ && $\mathbf{k}7_{53}$&$\kb{4^{4},\, 10^{3}}$&&& \\ \hline
\end{tabular} }

\end{document}